\documentclass[11pt,leqno]{article}
\usepackage[latin1]{inputenc}
\usepackage[T1]{fontenc}
\usepackage{amsmath}
\usepackage{amssymb}
\usepackage{amsfonts}
\usepackage{amsbsy}
\usepackage{amscd}
\usepackage{theorem}
\usepackage{epsf}
\usepackage{psfrag}
\usepackage{latexsym}
\allowdisplaybreaks

\def\build#1_#2^#3{\mathrel{\mathop{\kern 0pt#1}\limits_{#2}^{#3}}}
\def\noi{{\noindent}}
\def\be{\begin{equation}}
\def\ee{\end{equation}}
\def\ba{\begin{eqnarray*}}
\def\ea{\end{eqnarray*}}

\def\E{{\bf E}}

\def\cqfd{~~$\Box$\bigskip}

\def\o{\omega}

\def\ep{\varepsilon }

\def\O{\Omega}

\def\f{{\cal F}}

\def\n{{\cal N}}

\def\PP{{\bf P}}
\def\TT{{\bf T}}

\def\wh{\widehat}

\def\la{\longrightarrow}

\def\app{\frac{a_p}{p}}

\def\xbr{X^{br}}
\def\rxbr{\widehat{X}^{br} }
\def\hbr{H^{br} }
\def\mbr{M^{br} }
\def\wbrp{W^{br, p}}
\def\hbrp{H^{br , p}}
\def\mbrp{M^{br,p}}
\def\rwbrp{\widehat{W}^{br,p}}
\def\rlbrp{\widehat{L}^{br,p}}
\def\lbrp{L^{br,p}}
\def\xexc{ X^{exc} }
\def\hexc{H^{exc} }
\def\noi{\noindent}
\def\dem{{\sc Proof.} }

\def\card{{\rm Card\,}}

\def\build#1_#2^#3{\mathrel{\mathop{\kern 0pt#1}\limits_{#2}^{#3}}}

\newtheorem{theorem}{Theorem}[section]
\newtheorem{lemma}[theorem]{Lemma}
\newtheorem{proposition}[theorem]{Proposition}

{\theorembodyfont{\rmfamily}\newtheorem{remark}{Remark}[section]}
\newcommand{\R}{\mathbb{R}}

\newcommand{\Z}{\mathbb{Z}}

\newcommand{\Q}{\mathbb{Q}}

\newcommand{\D}{\mathbb{D}}

\newcommand{\NN}{\mathbb{N}}

\newcommand{\un}{\boldsymbol{1}}

\begin{document}

\title{A LIMIT THEOREM FOR THE CONTOUR PROCESS OF
CONDITIONED GALTON--WATSON TREES}
\author{by  Thomas Duquesne\\
{\it Universit\'e Paris 11, Math\'ematiques,
    91405 Orsay Cedex, France}}

\maketitle

\begin{abstract}

In this work, we study asymptotics of the genealogy of Galton--Watson processes conditioned on the
total progeny. We consider a fixed,
aperiodic and critical offspring distribution such that the rescaled Galton--Watson processes converges
to a continuous-state branching process (CSBP) with a stable branching mechanism of index
$\alpha  \in  (1, 2]$. We code the genealogy by two different processes: the contour process and
the height process that Le Gall and Le Jan recently introduced \cite{LGLJ1, LGLJ1}.
We show that the rescaled height
process of the corresponding Galton--Watson family tree, with one ancestor and conditioned on the total
progeny, converges in a functional sense, to a new process:
the normalized excursion of the continuous height process
associated with the $\alpha $-stable
CSBP. We deduce from this convergence an analogous limit theorem for the contour process.
In the Brownian case $\alpha =2$, the limiting process is the normalized Brownian excursion that codes
the continuum random tree: the result is due to Aldous who used a different method.

\vspace{1cm}

\noi {\it AMS 2000 subject classifications.} 60F17, 05G05, 60G52, 60G17.

\noi {\it Key words and phrases.} Stable continuous random tree, limit theorem,
conditioned Galton--Watson tree.

\vspace{1cm}

\end{abstract}

\section{Introduction.}

The analogues in continuous time of the Galton--Watson branching
processes (G-W processes) are the continuous-state branching processes (CSBP).
This class of Markov processes was originally introduced by Jirina and Lamperti (see \cite{Ji}
and \cite{La2}). These processes are
the only possible
weak limits that can be obtained from sequences of rescaled G-W processes (see \cite{La1}
or \cite{La3}). The properties of CSBP have been extensively studied
(see Grey \cite{Grey} or Bingham \cite{Bi2}). Lamperti has shown that a general CSBP can be
obtained from a L\'evy process without negative jump  by a random time change. The Laplace exponent $\psi$
of the L\'evy process is called the branching mechanism of the CSBP and it characterizes its law
via a differential equation solved by the Laplace exponent of the process.

 When one considers sequences of rescaled G-W processes with some fixed offspring distribution
 $\mu $ on $\NN $, the possible
limit processes are the CSBP with stable branching mechanism, that is, $\psi (\lambda ) =
c\lambda^{\alpha } $, for
some positive $c$ and $\alpha $ in $(0, 2]$ (see \cite{La3}).
In the case $\psi (\lambda ) =c \lambda^2 $,
the corresponding CSBP is the Feller diffusion.

  In this work we use some recent results concerning the genealogical structure of CSBP that can be found
in \cite{LGLJ1}, \cite{LGLJ2} and \cite{DuLG}. Our basic object is the G-W tree with offspring
distribution $\mu $
that can be seen as the underlying family tree of the corresponding G-W process started with one
ancestor; this random
tree is chosen to be rooted and ordered (see Neveu \cite{Ne} for a rigorous definition). If $\mu $
is critical or subcritical,
the G-W tree is almost surely finite and it can be coded by two different discrete processes: the
contour process and the height process
that are both defined at the beginning of Section 2.
These two processes are not Markovian in general but they can be written as a functional
of a certain left-continuous random walk whose jump distribution depends on $\mu $ in a simple way.

If a sequence of rescaled G-W processes converges to a CSBP with branching mechanism $\psi$, then
it has been shown in
\cite{DuLG}, Chapter 2, that the genealogical structure of the G-W processes
converges too. More precisely, the corresponding rescaled
sequences of contour processes and height processes, converge respectively to $(H_{t/2} )_{t\geq 0 }$
and $(H_t )_{t\geq 0}$, where the limit process $(H_t )_{t\geq 0}$ is the height process in continuous
time that has been introduced by Le Gall and Le Jan in \cite{LGLJ1}.
As in the discrete case, the height process is not Markovian in general but it can be written as a
functional of the L\'evy process without negative jump,
with Laplace exponent $\psi $, that plays the role of the left-continuous random walk.

  The case of a height process corresponding to a L\'evy process with finite variation paths is
  treated in \cite{LGLJ1}. It has
an interpretation in terms of queuing processes that has been used in some recent work of V. Limic
(see \cite{Li}). In the present
article,
we are only dealing with the case of the $\alpha $-stable branching mechanism, with $\alpha $
in $(1, 2]$. In that case the CSBP
is conservative and becomes extinct
almost surely. A general theorem implies that the corresponding height process is continuous
(see Theorem 4.7.
in \cite {LGLJ1} or \cite{DuLG}, Chapter 1) and the convergence of the rescaled discrete height processes
holds in a functional
sense. Furthermore, as explained in Section 3, the height process has a scaling property of
index $\alpha / (\alpha -1 )$. In the Brownian case $\alpha =2$, the height process is proportional
to a reflected standard
Brownian motion.

 In \cite{Al1} and \cite{Al2}, Aldous  introduced the continuum
random tree as the limit of rescaled G-W trees
conditioned on the total progeny, in the case where the offspring distribution
has finite variance. The continuum random tree is coded by a normalized Brownian excursion, in a
way similar to our
coding of discrete trees through the height process.
In the present work, we aim to extend Aldous' result to G-W trees with possibly
infinite variance offspring distribution. More precisely, we assume that the offspring distribution
of the G-W tree belongs
to the domain of attraction of a stable law with index $\alpha $ in $(1, 2]$. We then show that the
(suitably rescaled)
discrete height process of the G-W tree conditioned to have a large fixed progeny, converges in a
functional sense to the normalized excursion of the height process associated with the $\alpha $-stable
CSBP. This is the main result of the present work
and it is stated at the end of Section 3. We can think of our limiting process as the height process of
an infinite tree: by analogy, we call it
the $\alpha$-stable continuum random tree. In the case $\alpha = 2$, it coincides with
Aldous' continuum random tree. At the end of the Section 3, we also recall from \cite{DuLG}, Chapter 3,
the computation of finite dimensional marginals of the $\alpha$-stable continuum random tree.

 The last section is devoted to the proof of the limit theorem. Our approach
relies on an idea used by Kersting who introduced discrete bridges
in \cite{Ker} to study the convergence of rescaled G-W processes conditioned on the total
progeny, in the case of an infinite variance offspring distribution.
The limiting procedure is made easier in terms of discrete bridges thanks to their
good properties of absolute continuity with respect to the law of
the unconditioned random walk. In Section 4.1, we show that the height process of
the G-W tree conditioned on  its total progeny
has the same law as a certain functional of the discrete bridge. In the next section, we state a similar
result in the continuous setting. Then, we pass to the limit on functionals of discrete
bridges. The identification of the limit process as the normalized excursion
of the continuous height process involves several arguments that depend
on  continuity properties of the Vervaat transform (see \cite{Ve})
and on certain path-decompositions of the $\alpha$-stable L\'evy bridge that are due to Chaumont.

 \section{The coding of discrete Galton--Watson trees.}

In this section, we introduce the contour process and the
height process of a Galton--Watson
tree with a critical or subcritical offspring distribution. Each of these processes
provides a coding of the tree.
The height process can be written as a simple
functional of a left-continuous
random walk. This observation explains the definition of the continuous height process, that is given
in a forthcoming section. The results of this section are elementary and we refer
to \cite{LGLJ1} and \cite{DuLG} for details.

The trees considered in the present article are rooted ordered trees. Let us define them formally. We set
$ \NN^* = \{ 1, 2 , \ldots \}$ and
$$U=\bigcup_{n=0}^\infty (\NN^{*})^{n}  $$
where by convention $(\NN^{*})^{ 0}=\{\varnothing\}$. $U$ is the set of all possible words that can be
written with the elements of $\NN^* $. An element $u$ of $(\NN^* )^{ n}$ is written $u=u_1\ldots
u_n$, and we set $|u| =n$. If $u=u_1\ldots u_m$ and
$v=v_1\ldots v_n$ belong to $U$, we write $uv=u_1\ldots u_mv_1\ldots v_n$
for the concatenation of $u$ and $v$. In particular $u\varnothing=\varnothing u=u$.
We write $u<v$ for the lexicographical order on $U$: $\varnothing<1<11<12<121$
for example.

A rooted ordered tree $\tau$ is a subset of $U$ such that:

\ \ (i) $\varnothing\in \tau$.

\ (ii) If $v\in \tau$ and $v=uj$ for some
$j\in\NN^*$, then $u\in\tau$.

(iii) For every $u\in\tau$, there exists a number $k_u(\tau)\geq 0$
such that $uj\in\tau$ if and only if $1\leq j\leq k_u(\tau)$.

\noi We denote by $\TT$ the set of all trees. In the remainder, we see each vertex of a tree $\tau $
as an individual of some population
whose $\tau $ is the family tree and we shall often use a non-standard ``genealogical'' terminology
rather than the graph-theoretical one: for example, the individual $\varnothing $ is called
the ancestor of $\tau $.

Let us set some notation. Let $\tau_1 , \ldots , \tau_k $ be $k$ trees, the concatenation of
$\tau_1 , \ldots , \tau_k $, denoted by $[\tau_1 , \ldots , \tau_k ]$, is
defined in the following way: For $n \geq 1 $, $u=u_1u_2 \ldots u_n $
belongs to $[\tau_1 , \ldots , \tau_k ]$ if and only if $1\leq u_1 \leq k $
and $u_2 \ldots u_n $ belongs to $ \tau_{u_1 } $.

A leaf of the tree $\tau $ is an individual $u $ of $\tau $ that has no
child, as-to-say $ k_u (\tau ) =0 $. We denote
by ${\cal L}_{\tau } $ the set of all leaves of $\tau $.
If $\tau$ is a tree and $u\in \tau$, we define the shift of $\tau$ at $u$
by $T_u\tau=\{v\in U,uv\in\tau\}$.
Note that $T_u\tau\in \TT $. We denote by $\zeta (\tau ) =\card ( \tau ) $ the total progeny of $\tau $.
We write $u \preccurlyeq v $ if $v=uw$ for some $w$ in $U$
($\preccurlyeq $ is the ``genealogical'' order
on $\tau $). If $u\neq \varnothing $,
we use the notation $\overleftarrow{u }$ for the immediate
predecessor of $u$ with respect to $\preccurlyeq $, that can be seen
as the ``father'' of $u$ (thus $u=\overleftarrow{u }j $
for some positive integer $j$).
We also denote by $u\wedge v $ the youngest common ancestor of $u$ and $v$:
$$ u\wedge v = \sup \{   w\in \tau   :  w\preccurlyeq u \ {\rm and } \ w\preccurlyeq v \},$$
\noi where the supremum  is taken for the genealogical order.

 We now introduce the height process associated with a finite tree $\tau $. Let us denote by $u(0)=
 \varnothing < u(1) < u(2) <  \cdots
< u( \zeta (\tau ) -1 ) $ the individuals of $\tau $ listed in lexicographical order. The height process
$H(\tau )= (H_n ( \tau );0\leq n < \zeta (\tau ) ) $ is defined by
$$ H_n (\tau ) = |  u(n ) |    , \qquad 0\leq n < \zeta (\tau ) .$$
\noi The height process is thus the sequence of generations of the individuals of $\tau $ visited
in lexicographical order. It is easy to check that $H(\tau ) $ fully characterizes the tree.

  We also define the contour process associated with a tree $\tau $. We see $\tau $ embedded
in the oriented half-plane. We suppose
that the edges of $\tau $ have length one. Let us think of a particle visiting continuously each edge of
$\tau $ at speed one, from the left to the right: after
having reached $u(n )$, the particle goes to the individual
$u(n+1 )$, taking  the shortest way that consists first to move backward on the line of descent from
$ u(n ) $ to $ u(n ) \wedge  u( n+1 )$ and then, to move forward along the single
edge between $ u(n ) \wedge  u(n+1 ) $ to $u(n +1) $. The value $ C_t (\tau ) $ of the contour process
at time $t$ is the distance from the root to the position of the particle at time $t$. See Figure 1
for an example.

\vspace{3cm}

\begin{figure}[ht]
\psfrag{tau1}{$\tau$}
\psfrag{un}{$1$}
\psfrag{deux}{$2$}
\psfrag{un1}{$11$}
\psfrag{un2}{$12$}
\psfrag{un21}{$121$}
\psfrag{un22}{$122$}
\psfrag{vide}{$\varnothing$}
\psfrag{htau}{$H(\tau )$}
\psfrag{ctau}{$C (\tau )$}
\psfrag{0}{$0$}
\epsfxsize=12cm
\centerline{\epsfbox{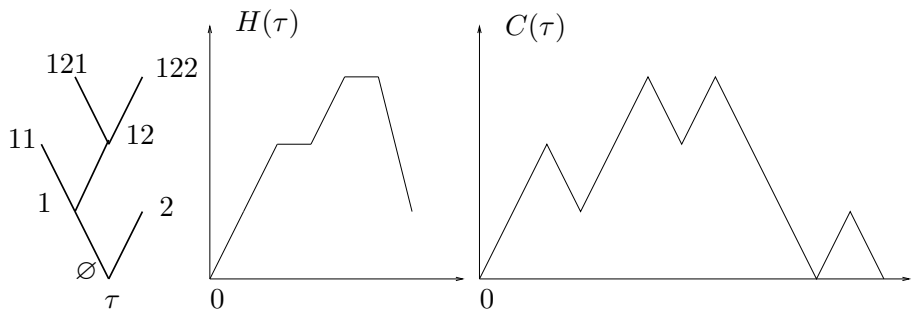}}
\caption{height process and contour process.}
\label{ramp.fig}
\end{figure}

More formally, we denote by $l_1 < l_2 <   \cdots   < l_p $ the $p $ leaves of $\tau $ listed
in lexicographical order. The contour
process $ ( C_t (\tau )   ;  t\in \R_+ ) $ is the piecewise linear continuous path with slope
equal to $ +1 $ or $ -1 $, that takes successive local extremes with
values: $ 0,  \mid l_1 \mid,\mid l_1 \wedge  l_2 \mid,\mid l_2 \mid,\ldots,
\mid l_{p-1 } \wedge  l_p \mid,\mid l_p \mid  $ and $0$.
Observe that the contour process visits each
edge of $\tau $ exactly two times. The contour process can be recovered from
the height process through the following transform. First set $b_n = 2n-H_n (\tau ) $, for $0\leq n <
\zeta (\tau )$ and $b_{\zeta (\tau ) }=
2(\zeta (\tau )-1)$. Then, observe that
$$ 0=b_0 < b_1 < \cdots <b_{\zeta (\tau ) -1 } < b_{\zeta (\tau )} = 2( \zeta (\tau )-1 ) .$$
For $n < \zeta (\tau ) -1 $ and
$t$ in $[b_n ,b_{n+1 }) $
\begin{equation}
\label{contour-height}
 C_t (\tau ) = \left\{
\begin{array}{ll}
H_n (\tau ) -(t-b_n ),&\qquad {\rm if \ }  t \in [b_n ,b_{n+1 }-1), \\
t -b_{n+1 } + H_{n+1 } (\tau ),&  \qquad {\rm if \ }  t\in [b_{n+1}-1 , b_{n+1 } ) ,
\end{array}  \right.
\end{equation}
and
$$ C_t (\tau ) = H_{\zeta (\tau )-1 } (\tau )-(t-b_{\zeta (\tau )-1}) \qquad {\rm if \ }
t\in [b_{\zeta (\tau )-1},
b_{ \zeta (\tau )}).$$

We can consider still another function coding
$\tau $, which is denoted by $(W_n(\tau);0\leq n <\zeta (\tau ) ) $ and defined by $ W_0 (\tau )=0 $ and
$$ W_{n+1 } (\tau ) -W_n (\tau ) = k_{u(n)}(\tau)-1,\qquad 0\leq n<\zeta(\tau ). $$
Observe that the jumps of $W(\tau ) $ are not smaller than $-1$.
The height process can be deduced from $W(\tau)$ by
the following formula (see Corollary 2.2. of \cite{LGLJ1}):
\begin{equation}
\label{height-walk}
H_n (\tau ) = \card \left\{ 0\leq j < n : W_j (\tau ) = \inf_{j\leq k\leq n}W_k(\tau)\right\},
\qquad 0\leq n< \zeta (\tau ).
\end{equation}
\noi As we will see in the next section, the
continuous height process is defined by analogy with this formula.

We now extend the definition of the height process, the contour process and the path $W$ to
a forest (i.e., a sequence) of finite trees: let $\varphi =
(\tau_p )_{p\geq 1 } $ be such a forest and set  $n_p = \zeta (\tau_1 ) +   \cdots   + \zeta (\tau _p ) $
with $ n_0 = 0 $. For any $p\geq 1 $, we define
$$ \left\{
\begin{array}{ll}
H_{n_p + k } (\varphi )   =   H_k (\tau_{p+1 } )   ,& \quad 0\leq k < \zeta( \tau_{p+1 } )   ,\\
W_{n_p + k } (\varphi )   =   W_k (\tau_{p+1 } ) -p  ,& \quad 0\leq k < \zeta( \tau_{p+1 } )   ,
\end{array} \right. $$
\noi and
$$ C_{t+2n_p -2p}(\varphi ) =
C_t (\tau_{p+1 } )   ,
\qquad  t \in [ 0, 2(\zeta (\tau_{p+1 } ) -1) ). $$
Observe that $\{ n_p ,   p\geq 0 \} $ is the set of integers $ k $ such that $ H_k (\varphi ) =0 $
or equivalently
$ W_k (\varphi ) < \inf_{0\leq j < k } W_j (\varphi )   $. Consequently, the excursions of
$ (H_n ( \varphi )   ;  n \geq 0 ) $
above $0$ (resp. the excursions of $( W_n (\varphi );\ n\geq 0)$ between
the successive times of decrease of its infimum) are the
$ (H_{n_p +k }  (\varphi )   ;  0\leq k < \zeta (\tau_{p+1 } ) ) $  (resp. the
$ (W_{n_p +k } (\varphi ) -p   ;  0\leq k < \zeta (\tau_{p+1 })  ) $ ). To the $p$th tree of
$\varphi $ corresponds the $p$th excursion of above level zero of $ H(\varphi ) $ and this
excursion coincides with
its height process.

\medskip

\begin{remark} In particular, this implies that (\ref{height-walk}) still
holds when $H(\tau )$ and $W(\tau )$ are replaced by $H(\varphi )$ , $W(\varphi ) $ respectively.
\end{remark}

\bigskip

 Let $\mu$ be a probability measure on
$\NN$. We assume that $\mu$ is
critical or subcritical:
$$\sum_{k=1}^\infty k\,\mu(k)\leq 1$$
and in order to avoid trivialities, we assume $\mu(1)<1$. The law of the Galton--Watson tree
with offspring distribution $\mu$
is the unique probability measure $P_\mu$ on
$\TT$ such that:

\ (i) $P_\mu(k_\varnothing=j)=\mu(j)$,\ $j\in \NN$.

(ii) For every $j\geq 1$ with $\mu(j)>0$, the shifted trees
$T_1\tau,\ldots,T_j\tau$
are independent under the conditional probability
$P_\mu(\cdot\mid k_\varnothing =j)$
and their conditional distribution is $P_\mu$.

Let $\varphi = ( \tau_p )_{p\geq 0 } $ be an i.i.d. sequence of G-W trees with offspring distribution
$\mu $. In general,
neither $H(\varphi ) $ nor $C(\varphi ) $ are Markovian. But it is easy to see that $W(\varphi )$ is a
random walk started
at zero; its jump distribution is $\nu (k) = \mu (k+1 )   ,   k\in \{ -1, 0, 1 , 2 \ldots \} $. This
property and
(\ref{height-walk}) imply the following proposition.
\begin{proposition}
\label{discrexc}
Let $\mu $ be a critical or subcritical offspring distribution. Let $(W_n   ;  n\geq 0 )$ be a random
walk started at $0$ with jump distribution $\nu (k ) = \mu ( k+1 ) $ , $ k\in \{ -1 , 0 , 1 \ldots \}$
defined under the probability measure $\PP $. Let us set
$ \zeta = \inf \{ n\geq 0   :   W_n =-1 \} $.
Define the
process $(H_n   ;  n\geq 0)$ by
$$ H_n = \card \left\{ 0\leq j <n   :  W_j =\inf_{j\leq l \leq n }W_l   \right\}.$$
\noi  Let $n\geq 1 $ be such that
$\PP (\zeta =n ) >0 $. The law of the
process $(H_n ;   0\le n <\zeta )$ under $\PP (. \mid  \zeta =n ) $,
is the same as the law of $H(\tau )$ under
$ P_{\mu } (  \cdot   \mid \zeta (\tau )=n )$.
\end{proposition}

\begin{remark}
The law of the G-W tree with a geometric offspring distribution conditioned to have its total progeny
equal to $n$, is the uniform probability measure on
the set of all ordered rooted
trees with $n$ vertices.
\end{remark}

\section{The $\alpha$-stable continuum random tree.}

\subsection{The height process.}

 We define the height process in continuous time by analogy with
(\ref{height-walk}). The role of the
left-continuous random walk is played by a stable
L\'evy process without
negative jump. In this section, we use several results about stable L\'evy processes and we refer
to \cite{Be}, Chapter VIII,  or to the original work of Chaumont \cite{Ch} and \cite{Ch97} for further
details.

  Let us denote by $ (\Omega , \f , \PP) $ the underlying probability space. Let $X$ be a process with
  paths in  $\D (\R_+ , \R ) $, the space of right-continuous with left limit
(c\`adl\`ag)
real-valued functions, endowed with the Skorokhod topology. We denote by $(\f_t )_{t\geq 0 } $ the
filtration generated by $X$ and augmented with the
$\PP $-null sets. We assume that $X$ is a stable L\'evy process
without negative jump with index $\alpha \in (1, 2]$. Then we have
$$\E \left[ \exp \left( -\lambda X_t \right) \right]   =   \exp \left( -c\lambda^{\alpha  } \right)  ,
\qquad \lambda >0   , $$
\noi for some positive constant $c$.
The process $(k^{-1/\alpha } X_{kt }   ;  t\geq 0 ) $ has the same law as
$(X_t   ;  t\geq 0 )$. Thanks to this scaling
property, we can take $c=1 $, without loss of generality in our purpose. When
$\alpha =2 $, the process $X$ is $1/\sqrt{2} $ times the standard Brownian motion on the line. When
$1<\alpha <2 $, the L\'evy measure of $X$ is
$$ \pi (dr )= \frac{\alpha (\alpha -1)}{\Gamma (2-\alpha ) } r^{-\alpha -1 } dr  . $$
\noi We use the following notation: for any $s<t $, we set
$$ I_{s,t } = \inf_{[s,t ]} X   ,\quad I_t = \inf_{[0, t ] } X   \quad {\rm and }
\quad S_t = \sup_{[0, t]} X   .$$
\noi Let us fix $t>0$. By analogy with the discrete
case, we want to define the height $H_t$
as the ``measure'' of the set
\begin{equation}
\label{zeroset}
\left\{ s\leq t ,  X_s=\inf_{s\leq r\leq t} X_r \right\}   .
\end{equation}
\noi To give a meaning to the word ``measure,'' we
use a time-reversal argument. Let $\wh X^{(t)}$ be the time-reversed process
$$ \left\{
\begin{array}{ll}
&\wh X^{(t)}_s=X_t-X_{(t-s)-},\qquad {\rm if}\ 0\leq s<t,\\
&\wh X^{(t)}_t=X_t.
\end{array} \right. $$
It is easy to check that
$(\wh X^{(t)}_s,0\leq s\leq t)
\overset{\rm(law)}{=} (X_s,0\leq s\leq t)$, that is refered to as the ``duality  property''.
We set $\wh S^{(t)}_s=\sup_{r\in [0,s]}\wh X^{(t)}_r$. Under the
transformation $s\la t-s$, the set (\ref{zeroset})
corresponds to
$$\left\{s\leq t,\,\wh X^{(t)}_s=\wh S^{(t)}_s\right\},$$
that is the zero set of the process $\wh S^{(t)}
-\wh X^{(t)}$ over $[0,t]$. Note that the
process $\wh S^{(t)} -\wh X^{(t)}$ has the same distribution as $S-X$.
However, the process $S-X$ is a Markov process. As
$\alpha \in (1, 2] $, the point $\{ 0\} $ is regular for itself with respect to this Markov process.
Hence, we can define the local time at $0$ of $S-X$ and denote it by $L=(L_t  ;  t\geq 0)$.
Note that $L$ is only defined up to a multiplicative
constant. Let us specify this normalization: Let $L^{-1}$ denote the right-continuous inverse of $L$,
$$L^{-1}(t)=\inf\{s>0  :  L_s>t \}.$$
\noi Both processes $(L^{-1}(t),t\geq 0)$
and $(S_{L^{-1}(t)},t\geq 0)$ are subordinators, called
respectively the ladder time process and the ladder
height process. The Laplace exponent of the ladder height process
is given by
$$ E\left[\exp \left( -\lambda S_{L^{-1}(t)} \right)\right]=\exp (-k {\lambda^{\alpha -1 }} ), $$
\noi where the positive constant $k$ depends on the normalization of $L$ (see \cite{Be}, Theorem VII-4). We
fix it by choosing $k=1 $.

\medskip

\begin{remark} Observe that in the Brownian case $\alpha =2 $, we have $L=S$.
\end{remark}
\medskip

   If $ 1< \alpha <2 $, we recall from \cite{LGLJ1} the following approximation of $L$. Let us
   denote by $(g_j, d_j )$ , $j\in J $ the
excursion intervals of $S-X$ above $0$. A classical argument of
fluctuation theory shows that the point measure
$$ \sum_{j\in J} \delta_{(L_{g_j} , \Delta S_{d_j} , \Delta X_{d_j} )}
(dldrdx ) $$
\noi is a Poisson measure with intensity $dl \pi (dx) \un_{[0,x]} (r) dr $
(see \cite{Ro84} or \cite{Be}, Chapter VI). Set
$$ \beta_{\varepsilon } = \int_{(\varepsilon , +\infty )} x\pi (dx)
= \frac{\alpha }{\Gamma (2-\alpha ) \ep^{\alpha-1 }}  . $$
By standard arguments we see that $\PP$-a.s. for every
$t\geq 0 $,
\begin{equation}
\label{L-approx1}
L_t = \lim_{\ep \rightarrow 0} \frac{1}{\beta_{\ep }} \card
\{ s\in [0, t]   :   S_{s-} < X_s   ;   \Delta X_s > \ep   \}   .
\end{equation}
\noi Thanks to this approximation, we can view
$L_t$ as a function of $(X_s ;  0\leq s \leq t )$.
Then we define the {\it height process in continuous time}, denoted by $H_t $, by the formula
$H_t= L_t ( \widehat{X}^{(t)} ) $.
In the Brownian case, the height process is $H_t = \widehat{S}^{(t)}_t = X_t -I_t $
and obviously has continuous paths.
If $ 1< \alpha <2 $, the
general theorem 4.7 of \cite{LGLJ1} shows that $H$ admits a continuous modification. Using the
 Fubini theorem, we deduce
from (\ref{L-approx1}) and from the duality property
that the limit
\begin{equation}
\label{H-approx1}
H_t = \lim_{\ep \rightarrow 0} \frac{1}{\beta_{\ep }} \card
\{ u\in [0, t ]   :   X_{u- } < I_{u,t }   ;  \Delta X_u > \ep   \}
\end{equation}
\noi holds $\PP$-a.s. on a set of values of $t$ of
full Lebesgue measure. We deduce from the
scaling property of $X$
and from the previous approximation formula that $H$ has a scaling property of index
$ \frac{\alpha }{\alpha -1 }$: For any $k>0 $
$$ \left( k^{ \frac{1}{\alpha} -1} H_{kt } ;
t\geq 0 \right)  \overset{\rm(law)}{=}
\left( H_t ;  t\geq 0 \right)  .$$

\subsection{ The normalized excursion of the height process.}

 Recall that $X-I$ is a strong Markov process and that $ 0  $ is regular for $X-I$. We may and will
 choose $-I$ as the local time of $X-I$ at level $0$. Let $(g_i , d_i ), i\in \cal{I}$ be the
excursion intervals of $X-I $ above $0$. Let us set
$$ \omega^{i }_s  = X_{g_i +s } -X_{g_i }   ,\qquad  0\leq s \leq \zeta_i = d_i - g_i   .$$
\noi The point measure
$$\n (dtd\omega ) = \sum_{i\in \cal{I} } \delta_{(-I_{g_i } , \omega^{i } ) }   $$
\noi is a Poisson measure with intensity $dt N(d\omega ) $. Here $N(d\omega )$ is a $\sigma$-finite
measure on the set of finite paths $(\omega (s); 0\leq s\leq
\zeta (\omega ) )$. Thanks to (\ref{H-approx1}), we see
that $H_t $ only depends on the excursion of $X-I $ straddling $t$. Thus we can use excursion theory
arguments in order to define the height process under the excursion measure $N$. We can also deduce from
(\ref{H-approx1}) that the excursions of $H$
above $0$ are almost surely equal to the $ H( \omega^{i }), i \in \cal{I}$ , with an evident
functional notation (see \cite{DuLG}, Chapter 1).

  We first have to define the normalized excursion of the $\alpha $-stable L\'evy process.
Let us simply denote by $\zeta =\zeta (\omega ) $ the lifetime of $\omega $ under $N (d\omega )$. A
standard result
of fluctuation theory says that $N( 1-e^{-\lambda \zeta } ) = \lambda^{1/\alpha } $ (see \cite{Be}).
Thus we have
$$ N( \zeta > t ) = \frac{1}{\Gamma (1-1/\alpha )} t^{  -\frac{1}{\alpha } }   .$$
Define for any $\lambda > 0$ the functional $S^{(\lambda )} $ by
$$ S^{(\lambda )} (\omega ) = \left( \lambda^{1/\alpha } \omega ( s/\lambda )   ;  0\leq s \leq
\lambda \zeta (\omega ) \right)   .$$
\noi Thanks to the scaling property of $X$, one can show that the image of $N(  \cdot \: \mid  \zeta >t ) $
under $ S^{(1/\zeta )} $ is the same for every $t>0 $. This law, defined on the
c\`adl\`ag paths with unit lifetime, is the law of the
{\it normalized excursion of $X$} denoted by $\PP^{exc} $. Informally $\PP^{exc } $ can be seen as
$N(\: \cdot \mid  \zeta =1 ) $
(see \cite{Be}, Chapter VIII). We
may assume that there exists a process $X^{exc }$ defined
on $(\Omega , \f , \PP )$ that takes values in $\D ( [0, 1 ], \R_+ ) $ and whose law under $\PP $
is $\PP^{exc}.$

   We recall Chaumont's path-construction of the normalized excursion of $X$ (see \cite{Ch}, \cite{Ch97}
or \cite{Be}, Chapter VIII): let $(\underline{g}_1 ,
\underline{d}_1) $ be the excursion interval of $X-I$ straddling $1$:
$$\left\{
\begin{array}{ll}
\underline{g}_1 &= \sup \{ s\leq 1   :  X_s = I_s \}   ,\\
\underline{d}_1 &= \inf \{ s>1   :  X_s = I_s   \}.
\end{array} \right.  $$
\noi We define $ \zeta_1 = \underline{d}_1 -\underline{g}_1 $, the length of
this excursion and we set
$$X^* = \left( \zeta_1^{-1/\alpha } \left( X_{\underline{g}_1 +
\zeta_1 s } - X_{\underline{g}_1 } \right)  ;  0\leq s\leq
1 \right)  . $$
Then, we have
\begin{equation}
\label{pathXexc}
X^{exc }   \overset{\rm(law)}{=}   X^*   .
\end{equation}

  Now, let us define the {\it normalized excursion of the height process}. In the Brownian case
  $\alpha =2 $, this
is the normalized excursion
of $X$. In the case $ 1<\alpha <2 $, the approximation
(\ref{H-approx1}) and the identity (\ref{pathXexc})
imply that the limit
$$ \zeta_1^{\frac{1}{\alpha } -1} H_{\underline{g}_1 +\zeta_1 t }
=\lim_{\ep \rightarrow 0}
\frac{1}{\beta_{\ep }}
\card \left\{ u\in [0, t]:X^{* }_{u- } < \inf_{[u, t ]} X^{* };\Delta X^{* }_u > \ep\right\} $$
\noi holds $\PP$-a.s. for a set of values of $t$ of full
Lebesgue measure on $[0, 1]$. So there exists a continuous process
$(H^{exc}_t   ;  0\leq t \leq 1  ) $ such that the limit
\begin{equation}
\label{Hexc-approx}
 H^{exc }_t =
\lim_{\ep \rightarrow 0}
\frac{1}{\beta_{\ep }}
\card \left\{ u\in [0, t]   :
X^{exc }_{u- } < \inf_{[u, t ]} X^{exc }   ;  \Delta X^{exc }_u >
\ep   \right\}
\end{equation}
holds $\PP$-a.s. for a set of values of $t$ of full Lebesgue
measure in $[0,1]$. The process $H^{exc } $ is
called the normalized excursion of the height
process. Moreover, we have
\begin{equation}
\label{pathHexc}
H^{exc }   \overset{\rm(law)}{=}   \left( \zeta_1^{\frac{1}{\alpha } -1}
H_{\underline{g}_1 +\zeta_1 s }   ;  0\leq s\leq 1 \right)   .
\end{equation}
\noi This result also holds in the Brownian case.

\subsection{The limit theorem.}

In this section, we state a limit theorem for
the rescaled discrete contour process and the rescaled discrete height process of a G-W tree
conditioned on its total progeny.
Before, we need to introduce some notation and to recall some results that are proved in \cite{DuLG},
 Chapter 2.

  Let $\mu $ be a critical or subcritical offspring distribution such that $ \mu (1 )<1 $ and let
  $(Z^p_n   ;  n\geq 0 )$ be a G-W process with offspring
distribution $\mu $, starting with $p $ ancestors: $Z^p_0 = p $. We let $\varphi = (\tau_p )_{p\geq 1 }$
be a sequence of i.i.d. G-W trees with offspring distribution $\mu $. By convenience, we denote by
$(H_n ;  n\geq 0 ) $, $(C_t ;  t\geq 0)$
and $(W_n ;  n\geq 0)$ the corresponding height process, contour process and random walk associated
with $\varphi $. As
was observerd in Section 2, $W$ is a left-continuous random walk with jump distribution $\nu $ defined by
$\nu (k)= \mu(k+1 ) $ , $ k\in \{-1, 0, 1 \ldots \} $.

  We assume that $\nu $ is in the domain of attraction of a stable law with index $\alpha \in (1, 2]$.
  The condition
$\nu ((-\infty , -1)) =0 $,
implies that the limit law is spectrally positive. Thus, there exists an increasing sequence of positive
real numbers
$(a_p )_{p\geq 0 }$ such that
$ a_p \rightarrow \infty $ and
\begin{equation*}
\label{attraction}
\tag*{({\bf H})}
\frac{1}{a_p } W_p \overset{(d)}{\longrightarrow } X_1
\end{equation*}
\noi where $X$ is a stable L\'evy process without negative jump with Laplace exponent
$\psi (\lambda )=\lambda^{\alpha }$ ,
$\alpha \in (1,2 ]$. Note that we have
automatically $a_p /p \rightarrow 0 $.
Grimvall has shown in \cite{Gr} that \ref{attraction} is equivalent to
\begin{equation}
\label{GWconv}
\left( \frac{1}{a_p} Z^{p}_{ [ \frac{p}{a_p } t ] }   ;  t\geq 0 \right)
\longrightarrow \left( Z_t   ;  t\geq 0 \right)   ,
\end{equation}
\noi where $(Z_t   ;  t\geq 0)$ is
a CSBP with branching mechanism $\psi (\lambda ) = \lambda^{\alpha } $.
Here and later,
{\it the convergence in distribution of processes always holds in the functional sense},
that is in the sense of the weak convergence of
the laws of the processes in the Skorokhod space
where they have their paths (which is meant by the symbol
$\overset{(d)}{\longrightarrow}$). We will use the notation
$\overset{(fd)}{\longrightarrow}$ to indicate weak
convergence of finite dimensional marginals.

Our starting points are Theorems 2.3.2 and 2.4.1 in \cite{DuLG}, that we
recall in our particular setting: under assumption \ref{attraction}, the following convergences hold:
\begin{equation}
\label{hauteurconve}
\left\{ \begin{array}{ll}
\left( \displaystyle\frac{a_p}{p} H_{[pt] }   ;  t\geq 0 \right)& \xrightarrow[\quad]{(d)}
\left( H_t   ;  t\geq 0 \right)   , \\[2ex]
\left( \displaystyle\frac{a_p}{p} C_{pt }   ;  t\geq 0 \right) & \xrightarrow[\quad]{(d)}
\left( H_{t/2 }   ;  t\geq 0 \right)   ,
\end{array} \right.
\end{equation}
\noi where $(H_t   ;  t\geq 0 )$ stands for the continuous height process associated with $X$.

 As in Section 2, we let $\tau $ be a G-W tree with offspring distribution $\mu $, under the
 probability measure $P_{\mu } $.
To simplify notation, we denote by $\zeta $ the total progeny of $\tau $. If $ \mu $ is assumed to
be aperiodic [i.e., $\rm{gcd} ( k\in \{ 1,2 , \dots \} : \mu (k) >0 ) =1 $], the conditional probability
$P_{\mu } (\cdot  \mid  \zeta=p )$
is well defined for $p \geq 1 $ sufficiently large.
Let $(H_n^{exc , p };  0\leq n\leq p ) $,
$(C_t^{exc , p };  0\leq t \leq 2p) $ and
$(W_n^{exc , p };  0\leq n \leq p )$ be
three processes defined on $(\Omega , \f , \PP )$
such that
$$\left( (H_n^{exc , p };  0\leq n< p ), (C_t^{exc , p }; 0\leq t < 2p-2), \left(W_n^{exc , p };  0\leq n <p
\right) \right) $$
has the same law as $( H(\tau ) , C(\tau ) , W(\tau ) )$
under $P_{\mu } (\cdot  \mid  \zeta =p )$ and such
that $ H^{exc, p }_p =0 $ , $C^{exc , p}_t =0 $ for
$t\in [2p-2, 2p ]$ and $W^{exc, p }_p =-1 $.
Let also $H^{exc} $ be the normalized excursion of the height process $H$ defined in the previous section.
The main goal of the present work is to prove the following limit theorem:
\begin{theorem}
\label{maintheo}
Assume $\ref{attraction}$ and that $\mu $ is aperiodic. Then, we have
$$ \left(
\frac{a_p}{p} H^{exc , p}_{[pt] }   ;  0\leq t\leq 1 \right)
\xrightarrow[\quad]{(d)} \left( H^{exc}_t   ;  0\leq t \leq 1  \right)  $$
\noi and
$$
\left(
\frac{a_p}{p}
C^{exc, p}_{pt }
  ;  0\leq t\leq 2\right) \xrightarrow[\quad]{(d)}
\left( H^{exc}_{t/2 }   ;  0\leq t\leq 2 \right)
  ,$$
\end{theorem}
\begin{remark}  Thanks to (\ref{contour-height}),
the second convergence of the theorem follows from the
first one: Set $b_n = 2n -H^{exc , p}_{n } $ ,
$0\leq n < p $ and $ b_p = 2p-2 $.
We deduce from (\ref{contour-height}) that, for
$0\leq n < p $,
\begin{equation}
\label{concontour1}
\sup_{b_n \leq t <  b_{n+1 } }\left| C^{exc, p }_t - H^{exc , p}_n \right| \leq
\left|  H^{exc, p}_{n+1 } - H^{exc, p}_n \right| +1   .
\end{equation}
\noi Define the random function $g_p  :  [0, 2p] \longrightarrow \NN $ by setting
$g_p(t) = n   $,  if   $t \in  [b_n , b_{n+1} )$ and $n<p $, and
$g_p (t) = p $, if $t\in  [2p-2 , 2p ] $. The definition
of $b_n $ implies
$$\sup_{0\leq t\leq 2p } \left| g_p(t) - \frac{t}{2} \right|  \leq \frac{1}{2} \sup_{ 0\leq k \leq p}
H^{exc, p}_{k} +1   .$$
\noi Set $f_p (t) = g_p ( pt) /p $. By (\ref{concontour1}), we have
$$ \sup_{t\in [0, 2]}\frac{a_p}{p} \left| C^{exc, p }_{pt } - H^{exc , p}_{p f_p (t) }
 \right| \leq \frac{a_p}{p} + \frac{a_p}{p}
\sup_{t\in [0, 1]} \left| H^{exc, p}_{[pt]+1 } - H^{exc, p}_{[pt]} \right|$$
\noi and
$$\sup_{t\in [0,2]} \left| f_p (t) - \frac{t}{2} \right|\leq \frac{1}{p} +
\frac{1}{2a_p} \sup_{ t\in [0,1]} \frac{a_p}{p} H^{exc, p}_{[pt]}   .$$
\noi Assuming that the first convergence of the theorem holds, we have
$$ \frac{a_p}{p} + \frac{a_p}{p}
\sup_{t\in [0, 1]} \left| H^{exc, p}_{[pt]+1 } - H^{exc, p}_{[pt]} \right |   \longrightarrow   0  $$
\noi and
$$ \frac{1}{p} +
\frac{1}{2a_p} \sup_{ t\in [0,1]}\frac{a_p}{p} H^{exc, p}_{[pt]}    \longrightarrow   0 $$
\noi in probability. Thus, the preceding bounds imply
$$\left(
\frac{a_p}{p}
C^{exc, p}_{pt }
  ;  0\leq t\leq 2\right) \xrightarrow[\quad]{(d)}
\left( H^{exc}_{t/2 }   ;  0\leq t\leq 2 \right)   . $$
\end{remark}

\medskip

\begin{remark} We denote by $Z^{exc , p}$
the G-W process started with one ancestor conditioned on having a total progeny equal to $p$.
 Under the same assumptions as
Theorem \ref{maintheo}, Kersting has proved in \cite{Ker} that
$(p^{-1} Z^{exc,p}_{[pt/a_p ]}   ;  t\geq 0 )$ converges in distribution in
$\D(\R_+, \R)$ to a process $Z^{exc } $ that is obtained from $X^{exc }$ by the Lamperti time
change. Theorem \ref{maintheo}
can be used to simplify Kersting's proof. More precisely, it implies Lemma 9 in \cite{Ker},
that is the key-argument showing
that the laws of $(p^{-1} Z^{exc,p}_{[pt/a_p ]}   ;  t\geq 0 )$
are tight.
\end{remark}

\begin{remark} If the offspring distribution has a finite variance,
then, $\alpha =2, H^{exc }$ is
proportional to the normalized Brownian excursion and Theorem \ref{maintheo} is
due to Aldous with a very different
proof (see \cite{Al2}). Let us mention that Bennies and Kersting proved a weaker
version of Aldous'theorem using
a method closer to our (see \cite{BK}).
\end{remark}

The convergence of Theorem \ref{maintheo} suggests that $H^{exc }$ is the height
process of a ``continuous tree.'' By
analogy with Aldous' continuum random tree, we call the limiting tree the
$\alpha$-stable continuum random tree that can be defined as a random compact metric
space in the following way: Each
$s\in [0,1]$ corresponds to a vertex at height $ H^{exc}_s $ in the $\alpha$-stable
continuum random tree. Let $t\in [0,1]$. The
distance in $\alpha$-stable continuum random tree between the two vertices corresponding
to $s$ and $t$ must be equal to
$$ d(s,t) =  H^{exc}_t + H^{exc}_s -2\inf_{u\in [ \min (s, t)  , \max (s, t)]}  H^{exc}_u   .$$
Then, we say the $s$ and $t$ are equivalent if and only if
$d(s, t)=0 $ and we denote it by $s\sim t $. We set
${\cal T } =  [ 0, 1 ] / \sim $ and we
define the $\alpha$-stable continuum random tree as the
(random) compact metric space $({\cal T } , d )$.
For a general theory, we refere to \cite{DMT96} and \cite{DT96}.

   For any $s\in [0,1]$ we denote by $\tilde{s}$ the corresponding vertex in $\cal{T}$;
    by analogy with the discrete case, we call
$\tilde{0}$ the root. The order on $\cal{T}$ induced by the order on $[0, 1]$ is the
continuous analogue of the lexicographical order
on discrete ordered trees. We can also define a "genealogical" order $\preceq $ on $\cal{T}$:
Let $\sigma , \sigma'\in \cal{T}$. Then we
say that
$$ \sigma \preceq \sigma' \qquad \rm{iff \ } d(\sigma , \sigma')= d(\tilde{0}, \sigma' )
-d( \tilde{0}, \sigma)   .$$
The set of leaves is the set of vertices that are maximal with respect to $\preceq$.
We denote it by $\cal{L}$. Here we give
some properties of $\cal{T}$ without proof (more general results are to be given in a forthcoming paper):

$\bullet $ $\PP$-a.s. the Lebesgue measure of $\{ s\in [0,1]   :  \tilde{s} \in \cal{L} \} $ is 1;

$\bullet $ $\PP$-a.s. the Hausdorff and packing dimensions of
$( {\cal T }, d )$ are both equal to $\frac{\alpha}{\alpha -1 }$;

$\bullet $  $\PP$-a.s. the Hausdorff and packing dimensions of $\cal{ T }\backslash \cal{L} $
are both equal to $1$.

Following Aldous \cite{Al1} and \cite{Al2}, we can define the finite dimensional marginals of $\cal{T}$.
Let us say a word about it: Aldous' first construction of the
(2-stable) continuum random tree was based on explicit formulas for the finite dimensional
marginals of this random tree.
Later, Aldous identified the continuum random tree as the tree coded by a
normalized Brownian excursion, in the sense of
\cite{Al2}. Le Gall \cite{LG2} provided a derivation of the
finite-dimensional marginals from properties of Brownian excursions.
A similar approach has been used in \cite{DuLG} to get the
finite-dimensional marginals of the $\alpha $-stable continuum random tree.
For sake of completeness let us explain how Theorem \ref{maintheo}
provides asymptotics for the finite-dimensional marginals of the G-W tree
conditioned on its total progeny.

Let $\tau $ be distributed under $P_{\mu }(\cdot \mid \zeta (\tau )=p )$ and fix
$k\leq p $. Let $(v_1, \ldots , , v_k )$ be a $k$-uple of distinct vertices of $\tau $. Aldous has
defined (Section 2 of \cite{Al2}) the $k$th
marginal of $\tau $ as the reduced subtree at $\{v_1 , \ldots , v_k  \}$
that is the (graph-theoretical) tree whose set of vertices $V$ is
$\{ v_i \wedge v_j ;  0\leq i \leq j \leq k \} \cup \{ \varnothing \} $ and whose
edges are all $(u, v)$, for $u$ and $v$ distinct in $V$ such that
$ u \preccurlyeq w\preccurlyeq v $ or $v \preccurlyeq w\preccurlyeq u $ occurs for $w \in V$
iff $w=u$ or $w=v $; furthermore, the length of the edge $(u, v)$ is
$\left| |u| -|v| \right| $.  Let us explain how the
$k$th marginal can be recovered from the height process of $\tau $.

  First we need to define what is a marked tree is: A marked tree is a
pair  $\theta = ( \tau , \{ h_v , v \in
\tau \} )$, where $\tau \in \TT $ and $ h_v \geq 0 $ for every $v \in \tau $. The number $h_v $
is interpreted as the lifetime of individual $v$ and $\tau $ is called the
skeleton of $\theta $. Let  $\theta_1 = ( \tau_1 , \{ h^1_v , v \in
\tau_1 \} ), \ldots , \theta_k = ( \tau_k , \{ h^k_v , v \in
\tau_k \} )$ be $k$ marked trees and $h\geq 0 $. The concatenation
of $ [\theta_1 , \ldots , \theta_k ]_h $ is the marked tree whose skeleton
is $[\tau_1 , \ldots , \tau_k ]$ and such that the lifetimes of vertices in
$\tau_i $ , $ 1\leq i \leq k $ become the lifetimes of the corresponding
vertices in $[\tau_1 , \ldots ,\tau_k ] $, and finally the lifetime of
$\varnothing $ in $[ \tau_1 , \ldots , \tau_k ]$ is $h$.

Assume that $k<p$  and let us explain how we deduce the $k$th marginals of
$\tau $ under $P_{\mu } (\cdot \mid \zeta (\tau )=p )$
from $H^{exc, p}$. Let
$\omega   :  [a,b] \rightarrow [0, +\infty )$ be a c\`adl\`ag function
defined on the subinterval $[a, b]$ of $[0, +\infty )$. Let $t_1 ,t_2 ,
\ldots , t_k \in [0 , +\infty ) $ be such that
$a\leq t_1 \leq  \cdots \leq t_k \leq b $. We first give the definition
of the marked tree associated to $\omega $ and $ t_1 , \ldots ,  t_k $.
For every $ a\leq u \leq v \leq b $, we set
$$ m(u,v)= \inf_{u\leq t\leq v } \o (t).$$
We will now construct a marked tree
$$ \theta (\o , t_1 ,\ldots t_k )= \left( \tau ( \o , t_1 , \ldots , t_k )
,  \{ h_v ( \o , t_1 , \ldots , t_k ) , v\in \tau \} \right) $$
\noi associated with the function $\o $ and the instants
$t_1 , \ldots , t_k $. We proceed by induction on $k$. If $k=1$,
$\tau ( \o , t_1 ) = \{ \varnothing \} $ and
$ h_{\varnothing } ( \o , t_1 ) = \o (t_1 ) $.

  Let $k\geq 2 $ and suppose that the tree has been constructed up to order
$k-1 $. Then there exists an integer $l \in \{ 1, \ldots , k-1 \} $ and
$l$ integers $1\leq i_1 < i_2 < \ldots < i_l \leq k-1 $ such that
$m(t_i , t_{i+1 } ) = m( t_1 , t_k ) $ if and only if
$ i\in \{ i_1 , i_2 , \ldots , i_l \} $. For every $q \in \{ 0 , 1 \ldots ,
l \} $, define $ \o^q $ by the formulas
$$ \begin{array}{lll}
\o^0 (t) = \o (t) - m( t_1 , t_k )   , \qquad t\in [a, t_{i_1}] , \\
\o^q (t) = \o (t) - m(t_1 , t_k )    , \qquad t\in [t_{i_q +1 } ,
t_{i_{q+1 } }] , \\
\o^l (t) = \o (t) - m(t_1 , t_k )   , \qquad t\in [ t_{i_l +1 } , b ] .
\end{array} $$
We then set
$$  \theta (\o , t_1 ,\ldots t_k )= \left[
 \theta (\o^0 , t_1 ,\ldots t_{i_1 }  ) ,  \theta (\o^1 , t_{i_1 +1} ,
\ldots t_{i_2 } ) , \ldots ,  \theta (\o^l , t_{i_l +1 } ,\ldots t_k )
\right]_{ m (t_1 , t_k)}.$$
This completes the construction of the tree by induction. Note that $l+1 $ is
the number of children of $\varnothing $ in
$\theta (\o , t_1 ,\ldots t_k ) $ and $ m(t_1 , t_k )$ is its lifetime. Figure 2 gives an example of a tree
$\theta ( \omega , t_1, \ldots , t_k ) $ when $k=4 $ and $[a, b]=[0,1]$.

\begin{figure}[ht]
\psfrag{s1}{$t_1$}
\psfrag{s2}{$t_3$}
\psfrag{s3}{$t_4$}
\psfrag{s4}{$t_2$}
\psfrag{zeta}{$1$}
\psfrag{te}{$t$}
\psfrag{hvide}{$h_{\varnothing }$}
\psfrag{hun}{$h_1$}
\psfrag{hdeux}{$h_2$}
\psfrag{zero}{$0$}
\psfrag{e(t)}{$\omega (t)$}
\epsfxsize=12cm
\centerline{\epsfbox{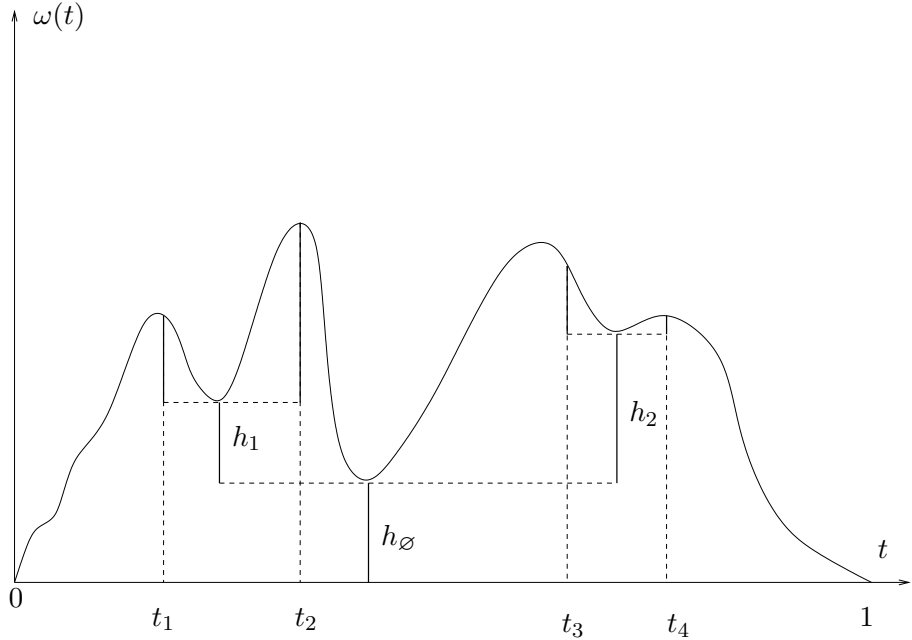}}
\caption{Reduced tree.}
\label{marginale}
\end{figure}

Let $ ( U^p_1 , U^p_2 , \ldots , U^p_k ) $ be independent of $ H^{exc , p } $ and
uniformly distributed on the set of all
$(n_1 , n_2 , \ldots , n_k )$ with $0\leq n_1 < n_2 < \cdots < n_k \leq p-1 $. From our
construction of the height process,
it should be clear that the $k$th marginal under $P_{\mu }  (\cdot |\zeta =p  ) $ is close to the tree
$\theta ( (H_{[pt]}^{exc , p })_{0\leq t\leq 1}   ,   U^p_1 , \ldots , U^p_k ) $
(in a sense that we do not make precise, but the reader can easily convince himself
that both trees have the same scaling limits when $p\rightarrow \infty $). On the other hand,
Theorem \ref{maintheo} implies that the rescaled trees $ \theta ( (\frac{a_p}
{p} H_{[pt]}^{exc , p })_{0\leq t \leq 1 }
,  \frac{U^p_1}{p} , \ldots , \frac{U^p_k }{p} ) $ converges in distribution to
$\theta ( H^{exc}, U_1, \ldots , U_k ) $, where
the $k$-tuple $(U_1 , \ldots , U_k )$ is independent of $ H^{exc} $ and distributed
according to the measure
$$ k! {\bf 1}_{ \{ 0 \leq u_1 < u_2 < \cdots < u_k \leq 1 \} } du_1 \cdots du_k.$$

We define $\theta ( H^{exc}, U_1, \ldots , U_k ) $ as the $k$th marginal of $\cal{T}$
The following theorem gives the law $\theta ( H^{exc}, U_1, \ldots , U_k ) $:

\begin{theorem}[Theorem 3.3.3 of \cite{DuLG}]\label{marginals}.  The law of
$\theta ( H^{exc} , U_1 , \ldots , U_k )$ is characterized
by the following properties$:$

\ {\rm(i)} The probability of a given skeleton $\tau
\in \{ \theta \in \TT :   |{\cal L}_{ \theta } | = k \quad and \quad
k_{u} (\theta ) \neq 1 ,   u\in \theta \} $ is
$$
\frac{k!}{\prod_{v\in \tau \backslash {\cal L}_{\tau } } k_v (\tau ) ! }
\frac{ \prod_{v\in \tau \backslash {\cal L}_{\tau } }
|(1-\alpha )(2-\alpha ) (3-\alpha ) \ldots ( k_v (\tau ) -1-\alpha )| }
{ |(\alpha -1 ) (2\alpha -1) \ldots ( (k-1) \alpha  -1 )| }   . $$

{\rm(ii)} Conditionally on the skeleton $\tau $, the marks
$(h_v )_{v\in \tau } $ have a density with respect to the Lebesgue measure
on $\R^{\tau }_{+} $ given by
$$  \frac{ \Gamma (k- \frac{1}{\alpha } ) }{ \Gamma ( \delta_{\tau } ) }
\alpha^{ \mid \tau \mid  }
\int_0^1 du   u^{ \delta_{\tau } -1 } q \left( \alpha \sum_{ v\in \tau }
 h_v   ,  1-u \right)   , $$

\noi where $ \delta_{\tau } = k-(1-\frac{1}{\alpha }) \mid \tau \mid  -
\frac{1}{\alpha } >0 $, and $q(s,u) $ is the continuous density at time
$s$ of the stable subordinator with index $ 1-\frac{1}{\alpha } $, that
is characterized by
$$ \int_0^{+ \infty } du e^{-\lambda s } q(s,u) = \exp \left( -s
\lambda^{1-\frac{1}{\alpha } }\right)   . $$
\end{theorem}

\begin{remark} In particular the skeleton of
$\theta ( H^{exc}, U_1, U_2, U_3 ) $ is equal to the discrete tree $\{ \varnothing ,
1,2,3 \}$ with probability $\frac{2-\alpha}{2\alpha -1}$.
Consequently, $\cal{T} $ has branching points of
order greater than $2$ if $\alpha <2$. General arguments (see \cite{DuLG}, Chapter 1)
imply that $\cal{T}$ has an infinite number of infinitely branching vertices.
\end{remark}

\section{Proof of the main theorem.}
The proof of Theorem \ref{maintheo} use Chaumont's result
on the Vervaat transform of the bridge of a $\alpha $-stable L\'evy process
(see \cite{Ch97} or \cite{Be}, Chapter VIII).
In this section, we explain how
the normalized excursion of the height process
is connected  (through the Vervaat transform) to
the height process associated with
the bridge of the L\'evy process. Before that, we need to establish some properties
in the discrete setting. This is the purpose of the following subsection.

\subsection{The discrete bridge.}

Let us start with some notation. We denote by
$\O_0 $ the set of all discrete-time finite paths in $\Z$:
$$ \O_0 = \bigcup_{n\geq 0 } \Z^{\{ 0, 1 ,
\ldots , n \} }   .$$
If $ w $ is in $\Z^{\{ 0 , 1 , \ldots , n \} } $ , we denote by $z(w) =n $ its lifetime.
Let $w$ be such that $z(w) \geq n $. We denote by respectively
$w^{(n)} $ and $\widehat{w}^n $, the shifted path and the time and space reversed path at time $n$:
$$ w^{(n)} (k) = w(k+n ) - w(n )   ,\qquad  0\leq k \leq z(w) -n$$
and
$$ \widehat{w}^n (k) = w(n ) - w( n-k)   ,\qquad   0\leq k\leq n.$$
\noi We set
$$ L_n (w ) = \card \left\{ 0 < j \leq n   :  w(j) =\sup_{0\leq k\leq j }
w(k)  \right\}   .$$
\noi We also define
$$ H_n (w) = L_n ( \widehat{w}^n ) = \card \left\{ 0 \leq  j < n   :  w(j) =\inf_{j\leq k\leq n }
w(k)  \right\}   . $$
\noi For any integer $a$, we define $ t (a, w )$ by
$$ t (a, w) = \inf \{ k \in [0, z(w) ]    :  w(k)
\geq a \}  $$
\noi (with the convention $\inf \varnothing =+\infty $). A careful counting
leads to the following formulas, valid
for any $0\leq m\leq z(w) -n $:
\begin{equation}
\label{combi1}
\left\{
\begin{array}{ll}
\displaystyle
H_{n+m} (w )   -\inf_{n\leq k \leq n+m } H_k (w)  = H_m \left( w^{(n)}\right), \\
\displaystyle
\inf_{n\leq k \leq n+m } H_k (w)  = L_n
(\widehat{w}^n ) -L_{ \beta (n,m ) } \left(\widehat{w}^n \right),
\end{array} \right.
\end{equation}
\noi where
\begin{equation*}
\beta (n, m) = \left\{
\begin{array}{lll}
\displaystyle
0, &{\rm if}  \displaystyle \inf_{0\leq k \leq m } w^{(n)} (k) \geq 0    , \\
\displaystyle
t \left( -\inf_{0\leq k \leq m } w^{(n)} (k)   ,
\widehat{w}^n \right) -1,   &{\rm if} \displaystyle
-\sup_{0\leq k \leq n }\widehat{w}^n (k) \leq \inf_{0\leq k \leq m } w^{(n)} (k) <0   , \\
\displaystyle
n &{\rm if} \displaystyle   \inf_{0\leq k \leq m } w^{(n)} (k) < -\sup_{0\leq k
\leq n }\widehat{w}^n (k)   .
\end{array}
\right.
\end{equation*}
Shortly written, we have
$$ \beta (n, m) = n\wedge \left( t \left(-\inf_{0\leq k \leq m } w^{(n)} (k)   ,
\widehat{w}^n \right)-1 \right)_+   .$$
We also set
$$ G( w)   =   \inf \left\{ 0\leq k\leq z(w):w(k) = \inf_{0\leq j\leq z(w)} w(j) \right\}.$$
\noi We now define the Vervaat transform $V_0
:   \O_0 \longrightarrow \O_0 $ by
$$ V_0( w ) (k) = \left\{
\begin{array}{ll}
w (k+G (w ) ) - \inf w,\quad &{\rm if } \ 0\leq k \leq z(w) - G(w)    , \\
w (k+G (w ) -z(w) ) + w (z(w)) - \inf w - w(0),\quad &{\rm if } \
 z(w) -G(w) \leq k \leq z(w)   .
\end{array}
\right. $$
\noi Observe that the path $V_0 (w) $ starts at $0$ and that its lifetime is
$z(w) $.

  Let us consider the random walk $W$ whose jump distribution is given
by $\nu (k) = \mu ( k+1 ) $ , $ k\in \{ -1 , 0 , 1 , \ldots , \} $. Recall
from Section 2 that $\zeta = \inf \{ k \geq 0   :  W_k  = -1 \} $. Set for
any positive integer $p$, $G_p = G( (W_k)_{ 0\leq k \leq p } )$.
A well-known result states that for any positive
integer $p$ ,  $p\PP ( \zeta =p )=\PP ( W_p =-1 ) $ (see \cite{Ot49}). In the remainder, we
assume that $\PP ( W_p =-1 )>0 $. Then,
$\PP ( \zeta =p ) > 0 $. We recall the classical result on random walk,
that is due to Vervaat (see \cite{Ve}):
\begin{equation}
\label{truevervaat}
V_0 (W) \ {\rm under} \ \PP ( \cdot \mid W_p =-1 )
\overset{\rm(law)}{=}   W \ {\rm under} \ \PP( \cdot \mid \zeta =p )   .
\end{equation}
This identity connects the discrete bridge of length $p$ with the
excursion conditioned to last $p$. We want to establish a similar identity
for the height process. To this end, we introduce the process
$M= (M_k )_{0\leq k\leq p } $ that is defined by the formula
$$ M_k = L_p\left(\widehat{W}^p\right) - L_{ \gamma_p (k) } \left( \widehat{W}^p \right),$$
where we have set $\gamma_p (k) = p\wedge ( t (-\inf_{0 \leq i \leq k}
W_i,\widehat{W}^p)-1)_{+} $.

Let us explain the intuition behind $M_k$: Consider
$\tau$ under $P(\cdot|\zeta=p)$. Let $\varnothing=u_0<u_1 < \ldots <u_{p-1}$
be the vertices of $\tau$ lexicographically ordered. Pick
$N$ at random in $\{0, 1, \ldots, p-1\}$ and assume that $N$ is independent
of $\tau$. We denote by $N(k)$ the integer of
$\{0, 1, \ldots, p-1\}$ equal to $N+k$ modulo $p$. Then $M_k$ has the same law
as $|u_{N}\wedge u_{N(k)}|$ that is the height of the common ancestor of $u_N$ and  $u_{N(k)}$.
We have the following proposition.

\begin{proposition}
\label{Vervaatdiscr} The law of the process $(V_0(W), V_0(H(W)+M))$ under
$\PP ( \cdot \mid W_p = -1 )$ is the same as that of $ (W, H(W) )$ under
$\PP ( \cdot \mid \zeta =p) $.
\end{proposition}
\dem
Set $W' = V_0  (W) $. Thanks to (\ref{truevervaat}), it is sufficient to
prove that $\PP ( \cdot \mid W_p =-1)$-a.s. $H(W' ) = V_0 ( H(W) +M ) $.
Let $0\leq l \leq G_p $. Applying (\ref{combi1}) with $w=(W_k' ;  0\leq k\leq p )$,
$n= p-G_p $ and $m= l$, we get
$$ H_{p-G_p +l } (W' ) = H_l \left( W^{'(p-G_p )}\right) +
L_{p-G_p } \left(\widehat{W'}^{ p-G_p }\right) - L_{ \beta' (p-G_p , l ) }
\left(\widehat{W'}^{ p-G_p }\right)$$
where $\beta' (p-G_p , l )= \beta  (p-G_p , l ) (W') $. However,
$$  W^{'(p-G_p )} = (W_k )_{0\leq k \leq G_p } \quad {\rm and} \quad
\widehat{W'}^{ p-G_p } = \left(\widehat{W}^p_k\right)_{0\leq k \leq p-G_p }   .$$
Then, $ H_l ( W^{'(p-G_p )} ) = H_l (W) $ and it is easily verified that
$L_p (\widehat{W}^p )=L_{p-G_p}  (\widehat{W'}^{ p-G_p })$
and $L_{\gamma_p (l)} (\widehat{W}^p ) =  L_{ \beta' (p-G_p , l ) } ( \widehat{W'}^{ p-G_p } ) $,
so that
$$ L_{p-G_p } \left(\widehat{W'}^{ p-G_p }\right) - L_{ \beta' (p-G_p , l ) }
\left(\widehat{W'}^{ p-G_p }\right) = M_l   .$$
So we have
\begin{equation}
\label{vervaat1}
 H_{p-G_p +l } (W' ) = M_l + H_l (W)   ,\qquad 0\leq l\leq G_p   .
\end{equation}

  Let us consider now $0\leq l \leq p - G_p $. We then have
$$ W'_l = W_{l + G_p } - W_{G_p } = W_{l+ G_p } -
\inf_{0\leq k \leq p } W_k   . $$
It easily follows that $H_l (W' ) = H_{l+ G_p } (W) $. But
$ M_{l+G_p } =0 $ because $\gamma_p (l+G_p ) = p $ [note that $ - \inf_{0\leq k \leq p }
W_k =  \sup_{0\leq k \leq p } \widehat{W}^p_k
+1 $,  $\PP (\cdot | W_p =-1)$-a.s.]. We conclude that
\begin{equation}
\label{vervaat2}
H_l (W' ) = H_{l+G_p } (W) + M_{l+G_p }   , \qquad 0\leq l
\leq p-G_p   .
\end{equation}

  Thanks to (\ref{vervaat1}) and (\ref{vervaat2}), we see that it only remains
to prove that $G (M+ H(W) ) = G_p $: First note that if $G_p =p $, we have $M_l =0 $ for every
$l\in[0, p]$ and $H(W)=V_0 \left( H(W) \right) $ in a trivial way. We can therefore suppose
$0<G_p<p $. Then it is easily seen that $\PP (\cdot |W_p =-1 )$-a.s. , for every $l\in [0, p]$,
$$ M_l = \card \left\{ 0\leq j<p   :  W_j =
\inf_{ j\leq k \leq p}  W_k \ {\rm and } \
W_j \leq -1 +\inf_{0\leq k \leq l } W_k \right\}   .$$
If $0\leq l < G_p $, then $ \inf_{0\leq k\leq l} W_k > \inf_{0\leq k\leq p } W_k $ and
thus $M_l >0 $ because we can
take $ j= G_p $ in the previous formula . On the other hand,  $H_{G_p } (W) + M_{G_p } =0 $ .
This proves $G ( M+ H(W)  ) = G_p $.~~$\Box$\bigskip

\subsection{Auxiliary processes.}
In this section we introduce the
L\'evy bridge $X^{br }$ that can be seen informally as the path $(X_t  ;  0\leq t\leq 1)$
conditioned to be at level zero
at time one. Standard arguments make this singular conditioning rigorous
and we refer to the original work of Chaumont
\cite{Ch}, \cite{Ch97} or to \cite{Be}, Chapter VIII, for the proofs. We also define the height process
associated with the bridge, denoted by $H^{br }$ and
the process $M^{br }$ that will play the role of
$M$ in continuous time.

  We denote by $p_t $ the continuous density of the law of $X_t $; it is characterized by
$$\int_{\R } \exp (-\lambda x ) p_t (x) dx = \exp (-t\lambda^{ \alpha } )   .$$
\noi For $0<t<1 $, the law of $(X^{br }_s   ;  0\leq s\leq t)$ is absolutely continuous with
respect to the law of
$(X_s   ;  0\leq s \leq t )$. More precisely, for any bounded continuous functional
$F$ defined on  $\D ([0,t], \R) $, we have
\begin{equation}
\label{abscont}
 \E \left[ F\left( X^{br }_s   ;  0\leq s \leq t \right)  \right] = \E \left[ F \left(
X_s   ;  0\leq s \leq t \right)  \frac{p_{1-t} (-X_t )}{p_1 (0) } \right]   .
\end{equation}
\noi It follows that
\begin{equation}
\label{rev}
\widehat{X}^{br }   \overset{\rm(law)}{=}   X^{br }
\end{equation}
where, for convenience, we denote by $\widehat{X}^{br}$
the process $X^{br }$ reversed at time $1$
( $\widehat{X}^{br }_t = -X^{br }_{1-t} ,\ t\in [0, 1] $).
Chaumont provides in \cite{Ch97} a path-construction for $X^{br }$: set $G=
\sup \{ t\in [0,1]  :   X_t =0 \} $,
the last passage time
at the origin on $[0,1]$ of the unconditioned process $X$. Let us set $\widetilde{X}= (G^{-1/\alpha }
X_{Gt})_{0\leq t \leq 1} $.
Chaumont has shown that
\begin{equation}
\label{pathbridge}
X^{br}  \overset{\rm(law)}{=}   \widetilde{X}
\end{equation}

  In the Brownian case $\alpha =2 $, we define the two processes
 $H^{br } $ and $L( X^{br })$ by setting
$$ H_t^{br }= X_t^{br}- I_t^{br } \quad {\rm }\quad
L_t( X^{br }) = S_t^{br}   ,\qquad t\in [0,1]   ,$$
with an evident notation for $I^{br } $ and $S^{br } $.

  If $1<\alpha <2 $, we
define $\hbr $ and $L(\xbr )$ by use
of the approximation formula (\ref{H-approx1}): By (\ref{abscont}) and a continuity
argument it is easy to check
that $\PP( \exists \delta >0   :  S^{br}_{1-\delta } = S_t^{br}
,   t\in [1-\delta , 1]) =1 $. Then, by (\ref{L-approx1}) and
(\ref{abscont}), it follows that we may define a continuous increasing process $L (X^{br }) $ by setting
$\PP$-a.s. for every $t$ in $[0,1]$,
\begin{equation}
\label{L-br-approx}
L_t (X^{br }) = \lim_{\ep \rightarrow 0}  \frac{1}{\beta_{\ep }} \card \left\{s\in [0,t]:
S^{br }_{s-} < X^{br }_{s};\Delta X^{br }_s >\ep \right\}   .
\end{equation}
  Next, by (\ref{H-approx1}), it follows that the limit
$$ G^{\frac{1}{\alpha} -1} H_{Gt}
= \lim_{\ep \longrightarrow 0}  \frac{1}{\beta_{\ep }}
\card \left\{ s\in [0,t]:\widetilde{X}_{s-} < \inf_{[s, t]} \widetilde{X};\Delta
\widetilde{X}_s >\ep \right\} $$
holds $\PP$-a.s. for a set of values of $t$ of full Lebesgue measure in
$[0, 1]$. Then, thanks to Chaumont's identity (\ref{pathbridge}) we
can show that there exists a continuous process $(H_t^{br};  0\leq t
\leq 1 )$ such that the limit
\begin{equation}
\label{H-br-approx}
H_t^{br} =  \lim_{\ep \longrightarrow 0}  \frac{1}{\beta_{\ep }} \card \left\{ s\in [0,t]:
X^{br }_{s-} < \inf_{[s,t]} X^{br }   ;  \Delta X^{br }_s >\ep   \right\}
\end{equation}
holds $\PP$-a.s. for a set of values of $t$ of full Lebesgue measure in
$[0, 1]$. We also have
\begin{equation}
\label{Hpathbridge}
H^{br}  \overset{\rm(law)}{=}   \left( G^{\frac{1}{\alpha} -1 }
H_{Gt} \right)_{0\leq t \leq 1}   .
\end{equation}
And by (\ref{abscont}), it follows that, for $0<t<1$,
\begin{equation}
\label{Habscont}
\E \left[ F\left( H^{br }_s   ;  0\leq s \leq t \right)  \right] = \E \left[ F \left(
H_s   ;  0\leq s \leq t \right)  \frac{p_{1-t} (-X_t )}{p_1 (0) }
\right]   .
\end{equation}

\bigskip

\begin{remark} Equations (\ref{Habscont}) and (\ref{Hpathbridge}) both
hold in the Brownian case.
\end{remark}

\bigskip

  We now define the Vervaat transform in continuous time, denoted by
$V: \D ([0,1] , \R ) \longrightarrow  \D ([0,1] , \R )$:
For any $\omega $ in $\D ([0,1], \R )$, we set
$g_1 (\omega ) = \inf \{ t\in [0,1]   :  \o (t-) \wedge
\omega (t) =\inf_{[0,1]} \omega \} $. Then, we define $V$ by
$$ V(\omega ) (t) = \left\{
\begin{array}{ll}
\omega (t+g_1 (\omega ) ) - \inf_{[0,1]} \omega,\quad &{\rm if } \ t+ g_1 (\omega ) \leq 1   , \\
\omega (t+g_1 (\omega ) -1 ) + \omega (1) - \inf_{[0,1]} \omega  -\o (0),\quad
&{\rm if } \ t+ g_1 (\omega ) \geq 1   .
\end{array} \right. $$
  Thanks to (\ref{abscont}), it is easy to see that the bridge $X^{br }$ reaches its
  infimum almost surely at a unique random time (that must be
$g_1 (X^{br })$ and that is uniformly distributed in $[0,1]$).
The bridge is connected to the normalized excursion $X^{exc}$ through the Vervaat transform
\begin{equation}
\label{vervaatbridge}
V(X^{br } )   \overset{\rm(law)}{=}  X^{exc }  .
\end{equation}
(For a proof, see Chaumont \cite{Ch97} or Bertoin \cite{Be}, Chapter VIII.)
Next, we define the analogue of $M$ in continuous time:
For any $\omega $ in $\D( [0,1], \R)$ and any positive real number $x$,
let us denote by $T_x (\omega )$, the first passage time above $x$:
$$ T_x (\omega )= \inf \{ t\geq 0   :  \omega (t)
\geq x \}   , $$
\noi (with the convention:  $\inf \varnothing = +\infty $ ). For any
$1< \alpha \leq 2 $ , $  L_t (\widehat{X}^{br } )$ is well defined thanks to (\ref{rev})
and (\ref{L-br-approx}). So we can set
$$ B_x = L_1 (\widehat{X}^{br } )- L_{1\wedge  T_x (\widehat{X}^{br } ) } ( \rxbr )   ,
\quad  x\geq 0  $$
\noi and we define $\mbr $ by
$$ \mbr_t = B_{-\inf_{[0,t]} \xbr }   ,  0\leq t \leq 1   .$$
\noi The following proposition is an analogue in continuous time of Proposition \ref{Vervaatdiscr}.

\begin{proposition}
\label{vervaat}
The processes  $\hbr $,  $\mbr $ and $B$ have the following properties$:$

\ \ {\rm (i)} $\PP$-a.s. $(B_x  ;  x\geq 0)$ is a nonnegative and nonincreasing continuous process.
Furthermore
we have $B_x =0 $ if and only if $   x  \geq -\inf_{[0,1]} \xbr.$

\ {\rm(ii)} $\PP$-a.s. $( \mbr_t + \hbr_t   ;  0\leq t\leq 1 )$ is a nonnegative continuous process
that attains its minimal value $0$ at a unique instant.

{\rm(iii)}   $V( \mbr +\hbr )  \overset{\rm(law)}{=}   \hexc.$
\end{proposition}

\dem
Thanks to Chaumont's result (\ref{vervaatbridge}), we can  assume that $\xexc $
and $\xbr $ are related in the following way:
\begin{equation}
\label{presquesur}
\left\{
\begin{array}{ll}
\xexc_s = \xbr_{g_1 +s } -\xbr_{g_1 }   ,&\quad  0\leq s \leq 1-g_1   ,\\
\xexc_s = \xbr_{s+g_1 -1 } - \xbr_{g_1 }   ,&\quad  1-g_1 \leq s \leq 1   ,
\end{array}  \right.
\end{equation}
\noi where we have set $ g_1 = g_1 (\xbr )$.

  In the Brownian case $\alpha =2 $, we have $ L_t (\widehat{X}^{br } )= \sup_{s\leq t}
  \widehat{X}^{br }_s  $. It easily follows that
$B_x = (-x- I^{br }_1)_+ $, $ M^{br }_t =  I^{br }_t  - I^{br }_1 $ and $ M^{br }_t +
H^{br }_t = X^{br }_t  - I^{br }_1 $,
for $x\geq 0 $ and $ 0\leq t\leq 1 $. Assertions (i) and (ii) follow immediately,
and (iii) is a direct consequence of
(\ref{presquesur}).


   From now on, we assume that $ 1<\alpha <2 $. Let us prove
(i) first. Recall that $(S_{L^{-1 }_t}; t\geq 0)$ is a stable subordinator with index
$ \alpha -1 $. Hence, its right-continuous inverse $(L_{T_x (X)}   ;  x\geq 0)$ is
$\PP$-a.s. continuous. If $x\geq S_1 $, then $T_x (X) \geq 1 $ and $L_1 =L_{1\wedge  T_x (X) } $. However,
for any positive rational $q$, the Markov
property for $X$ implies that $T_q (X) $ is an increase time for $L $. Then
$$ L_{T_q (X)} < L_1  \qquad {\rm on } \ \{ q<S_1 \}.$$
\noi Hence, $(L_1 -L_{1\wedge  T_x (X)};x\geq 0 )$ is $\PP$-a.s. a nonincreasing
and nonnegative continuous process that vanishes if and only if $x \geq S_1$.

Let $t<1$. We can use property (\ref{abscont}) to show that $\PP$-a.s. the process
$(L_{t\wedge  T_x (\xbr )} (\xbr )   ;  x\geq 0 ) $ is continuous and $ L_{T_x (\xbr )} (\xbr ) < L_t (\xbr )$
if and only if  $x< \sup_{[0,t]} \xbr $.
But $\PP$-a.s. there exists $t\in [0,1) $ such that $\sup_{[0,t]} \xbr =\sup_{[0,1]} \xbr $ and so
$$ L_{t\wedge  T_x (\xbr ) } (\xbr ) = L_{1\wedge  T_x (\xbr )}(\xbr )   ,\qquad   x\geq 0   .$$
Hence, we have proved that $\PP$-a.s. the process $( L_{1\wedge  T_x (\xbr )}  (\xbr )   ;  x\geq 0 ) $
is continuous
and $ L_{T_x (\xbr )} (\xbr ) < L_1 (\xbr ) $ if and only if  $x< \sup_{[0,1]} \xbr $. Then, (i)
follows from the duality property (\ref{rev}). Then, the  continuity
of $M^{br }$ follows from the continuity of $I^{br}$.

Recall that $X^{exc}$ and $X^{br }$ are related by (\ref{presquesur}). We now establish the a.s. identity
\begin{equation}
\label{vervaatpresq}
H^{exc} = V( M^{br }+H^{br})   .
\end{equation}
First, observe that if $ t> g_1 $, then the conditions
$\xbr_{s-} < \inf_{[s,t]} \xbr $ and $s\in [0, t] $
imply that $ s\geq g_1 $. Thanks to the approximations (\ref{Hexc-approx})
and (\ref{H-br-approx}), and the
continuity of the processes $H^{exc}$ and $H^{br}$, we easily verify
that, $\PP$-a.s.,
\begin{equation}
\label{vervaatdeux}
\hexc_t = \hbr_{g_1 +t }, \qquad 0 \leq t \leq  1-g_1.
\end{equation}
However (i) and the definition of $M^{br} $ imply that
$\PP$-a.s. $M^{br}_t =0 $ for any $t$ in $[g_1 , 1]$. Then, by (\ref{vervaatdeux}), it follows that
\begin{equation}
\label{verv}
\hexc_t = \hbr_{g_1 +t  } +\mbr_{g_1 +t}   ,
\qquad 0\leq t \leq 1-g_1    .
\end{equation}

Next we have to prove
\begin{equation}
\label{attention}
\hexc_t = \hbr_{g_1 +t -1 } +\mbr_{g_1 +t-1}, \qquad 1-g_1 \leq t \leq 1.
\end{equation}
Set, for any $t> 1-g_1 $,
$$ \gamma_t = \sup \left\{ s<1:\xbr_s \leq I^{br}_{t+g_1 -1} \right\} = 1-
T_{-I^{br}_{t+g_1 -1 } } (\rxbr).$$
We also define, for any $ 0\leq s'\leq t' \leq 1 $
$$ N^{br}_{\ep }(s',t') = \card \left\{ u\in [0, s']:
X^{br }_{u-} < \inf_{[u,t']} X^{br }   ;  \Delta X^{br }_u >\ep\right\}   $$
If $1-g_1 < t  $, observe that
$$ \card \left\{ s\in [0, t]:X^{exc }_{s-} < \inf_{[s,t]} X^{exc };\Delta X^{exc }_s >\ep\right\}
= C_1 + C_2 + C_3,$$
where
\begin{equation*}
\left\{
\begin{array}{ll}
\displaystyle
C_1 &=  \displaystyle \card \left\{ s\in [0, \gamma_t - g_1 ]   :
X^{exc }_{s-} < \inf_{[s,t]} X^{exc }   ;  \Delta X^{exc }_s >\ep   \right\} =  N^{br}_{\ep }
(\gamma_t , 1)   , \\[2ex]
\displaystyle
C_2 &= \displaystyle \card \left\{ s\in (\gamma_t - g_1 , 1-g_1 )   :
X^{exc }_{s-} < \inf_{[s,t]} X^{exc }   ;  \Delta X^{exc }_s >\ep   \right\} = 0   ,\\[2ex]
\displaystyle
C_3 &=  \displaystyle \card \left\{ s\in [ 1-g_1 , t]   :
X^{exc }_{s-} < \inf_{[s,t]} X^{exc }   ;  \Delta X^{exc }_s >\ep   \right\} \\[2ex]
\displaystyle
&= \displaystyle N^{br}_{\ep } (t+g_1-1  , t+g_1 -1 )   .
\end{array} \right.
\end{equation*}
Thus,
\begin{equation}\begin{split}
\label{vervaat3}
&\card \left\{ s\in [0, t]   :
X^{exc }_{s-} < \inf_{[s,t]} X^{exc }   ;  \Delta X^{exc }_s >\ep   \right\}\\
&\qquad =  N^{br}_{\ep } (\gamma_t , 1) +
N^{br}_{\ep } (t+g_1-1  , t+g_1 -1 )   .
\end{split}\end{equation}
By approximation formula (\ref{L-br-approx}) applied to $L(\rxbr )$
it follows that $\PP$-a.s. for any $ 0 \leq t\leq 1$,
$$\lim_{\ep \rightarrow 0} \frac{1}{\beta_{\ep }} N^{br}_{\ep } (\gamma_t , 1 ) = L_1
(\rxbr ) -L_{1-\gamma_t }
(\rxbr )   = \mbr_{t+g_1 -1 }   .$$
But  approximation (\ref{Hexc-approx}) of $H^{exc } $ and  approximation
(\ref{H-br-approx}) of $H^{br}$ imply that the limits
$$  H^{exc}_t =
\lim_{\ep \rightarrow 0} \frac{1}{\beta_{\ep }}
\card \{ s\in [0, t]   :
X^{exc }_{s-} < \inf_{[s,t]} X^{exc }   ;  \Delta X^{exc }_s >\ep   \} $$
$$ H^{br}_{t+g_1 -1 } =
\lim_{\ep \rightarrow 0} \frac{1}{\beta_{\ep }}
N^{br}_{\ep } (t+g_1-1  , t+g_1 -1 ) $$
hold $\PP$-a.s. on a set of values of $t$ of full Lebesgue measure in
$[1-g_1 , 1]$. Then, (\ref{attention}) follows from (\ref{vervaat3})
and the continuity of
$M^{br }$, $H^{exc }$ and $H^{br} $.

 It remains to shows that $\mbr +\hbr $ reaches its infimum at the unique time $g_1 $. If $s> g_1 $, then
$\mbr_s =0 $ by (i) and $\hbr_s = \hexc_{s-g_1 } >0 $. If $g_1 >s$, then
$$ \inf_{[0,s]} \xbr   >   \xbr_{g_1 } =  \inf_{[0,1]} \xbr $$
\noi and (i) implies that $\mbr_s >0$. Finally, (ii) follows from the obvious fact
$ \mbr_{g_1 }=\hbr_{g_1 } =0 $.~~$\Box$\bigskip

  We now explain how the auxiliary processes $\xbr $, $\hbr $ and $\mbr $
are used in the proof of Theorem \ref{maintheo}:
Let $(W^{br,p}_t; t\in [0,1] )$ be a process whose distribution is the law of
$(1 /a_p  W_{[pt ]} ; 0\leq t\leq 1) $ under $\PP( \cdot \mid W_p =-1 )$.
Simultaneously with  $W^{br,p}$ we can introduce the processes $H^{br,p},
M^{br,p} , \widehat{W}^{br , p} , L^{br , p}  $ and $\widehat{L}^{br , p} $
(which can be written as functionals of $W^{br , p} $) that are such that
$$ \left( W^{br, p }, H^{br,p},
M^{br,p} , \widehat{W}^{br , p} , L^{br , p}, \widehat{L}^{br , p} \right)$$
has the same law as
$$ \left( \frac{1}{a_p} W_{[pt ]} ,
\frac{a_p}{p} H_{[pt ]} ( W ) ,
\frac{a_p}{p} M_{[pt]} ,
\frac{1}{a_p} \widehat{W}_{[pt ]} ,
\frac{a_p}{p} L_{[pt]} (W),
\frac{a_p}{p} L_{[p t]} (\widehat{W}^p )
\right)_{0\leq t\leq 1 } $$
under  $\PP( \cdot \mid W_p =-1 )$.
We have the following proposition.
\begin{proposition}
\label{propaux}
Under the assumptions of Theorem $\ref{maintheo},$ we have
$$ \left( \wbrp , \hbrp , \mbrp \right) \xrightarrow[p\rightarrow \infty]{(d)}
\left( \xbr , \hbr ,\mbr \right)  $$
\end{proposition}

\bigskip

Let us complete the proof of Theorem \ref{maintheo} thanks to Proposition \ref{propaux}
whose proof is postponed to the next section.

\bigskip

{\sc Proof Theorem \ref{maintheo}.}
First, it easy to deduce
from Proposition \ref{discrexc}, from the definition of the discrete Vervaat transform and
from Proposition
\ref{Vervaatdiscr} that
\begin{equation}
\label{continu}
\left(
V\left( W^{br,p} \right) ,  V\left( M^{br,p} + H^{br, p} \right) \right)
  \overset{\rm(law)}{=}
\left( \left(\frac{1}{a_p} W^{exc, p}_{[pt]}\right)_{0\leq t\leq 1}   ,
\left(\frac{a_p}{p} H^{exc,p}_{[pt]}\right)_{0\leq t\leq 1} \right)
\end{equation}
   Then, we need to prove some continuity property of $V$:
Let $(\omega_n )_{n\geq 0 }$ be a sequence of
paths in $\D( [0,1] , \R )$ that converges to $\omega $ for the Skorokhod topology.
If $\omega $ in continuous, then the convergence holds uniformly on $[0,1]$ :
$$\lim_{n\rightarrow \infty } \sup_{0\leq t\leq 1} \mid \omega_n (t) - \omega (t) \mid  =0 $$
\noi (see Jacod and Shiryaev \cite{Ja},
Chapter VI). Then, if we assume furthermore that $\omega $ attains its minimum at a
unique instant, it is easily seen that   $\lim g_1 (\omega_n ) = g_1 (\omega ) $. Thus,

\begin{equation}
\label{conti}
\lim_{n \rightarrow +\infty } \sup_{t\in [0,1]} \mid  V(\omega_n ) (t) - V(\omega )(t) \mid  = 0   .
\end{equation}
This shows that $V$ is continuous at any continuous path $\omega $ in $\D( [0,1] , \R ) $ that
attains its minimum at a unique
time. This observation combined with Proposition
\ref{vervaat}(ii) and Proposition \ref{propaux}, shows that
\begin{equation}
\label{deuxieme}
V\left( M^{br,p} + H^{br, p} \right)
\xrightarrow[\quad]{(d)}
V\left( M^{br} + H^{br} \right)   .
\end{equation}
Then Theorem \ref{maintheo} follows from (\ref{continu}) and from
Proposition \ref{vervaat} (iii).\cqfd

\subsection{Proof of Proposition \ref{propaux}. }

We first prove the following lemma for unconditioned processes.
\begin{lemma}
\label{jointconv}
Under the assumptions of Theorem $\ref{maintheo},$ the following joint convergence holds:
$$ \left( \frac{1}{a_p} W_{[pt]} , \frac{a_p}{p} L_{[pt]} (W), \frac{a_p}{p} H_{[pt]} (W) \right)
\xrightarrow[p\rightarrow \infty]{(d)}
\left( X, L, H \right) $$
\end{lemma}

\dem
A classical result on random walks shows that assumption \ref{attraction}
implies the convergence of $(1/a_p W_{[pt]}   ;  t\geq 0 )$ to $X$ in distribution
in $\D(\R_+ , \R)$ (see Jacod and Shiryaev \cite{Ja},
Chapter VII). Theorem 2.3.2 in \cite{DuLG} [recalled in (\ref{hauteurconve})] shows
that the rescaled height process $ ( \frac{a_{p}}{p} H_{[pt]} (W)   ;  t\geq 0 ) $ converges to $H$
in distribution in $\D (\R_+ ,\R)$
under assumption \ref{attraction}.

As a first step towards the proof of the convergence of rescaled process, it is also proved
in \cite{DuLG} (see Theorem 2.2.1) that
\begin{equation}
\label{finite}
 \left( \frac{a_{p}}{p} L_{[pt]} (W)\right)_{t\geq 0 }
\xrightarrow[p\rightarrow \infty]{(fd)} (L_t)_{t\geq 0 }   .
\end{equation}
As $L$ is a continuous nondecreasing process, a standard argument show
that the convergence (\ref{finite}) actually holds in distribution in
$\D(\R_+ , \R )$. Thus, the laws of the processes
$$\left( \frac{1}{a_p} W_{[pt]} , \frac{a_p}{p} L_{[pt]} (W), \frac{a_p}{p} H_{[pt]}
(W)   ;  t\geq 0 \right) $$
\noi are tight in the space of probability measures on $\D( \R_+ , \R^3 )$.

  If we look carefully at Theorem 2.2.1 in \cite{DuLG}, we see that the proof actually
  gives a stronger result than the weak convergence
of the finite dimensional marginals of the rescaled height process:
By the Skorokhod representation theorem, we can find a sequence of random walks $W^{'p}$, $p\geq 1 $,
each with the same law as $W$, and a L\'evy process $X'$ with
$$\left( \frac{1}{a_p} W^{'p}_{[pt]}   ;  t\geq 0 \right)
\longrightarrow (X'_t   ;  t\geq 0), $$
$\PP$-a.s. for the Skorokhod topology. Then the proof of Theorem 2.2.1 in \cite{DuLG} shows that
\begin{equation}
\label{convproba}
\frac{a_p}{p} L_{[pt]}(W^{'p} ) \longrightarrow L_t (X') \qquad {\rm and } \qquad \frac{a_p}{p}
H_{[pt]}(W^{'p} )\longrightarrow H_t(X')
\end{equation}
\noi in probability for every $t\geq 0$ (with an evident notation for $L(X') $ and $H(X')$).
It follows that
the only possible weak limit for the laws of
$$\left( \frac{1}{a_p} W_{[pt]} , \frac{a_p}{p} L_{[pt]} (W), \frac{a_p}{p} H_{[pt]}(W);t\geq 0\right) $$
is that of $(X, L , H )$ and the lemma is proved.\cqfd

\begin{lemma}
\label{lemme1}
Under the assumptions of Theorem $\ref{maintheo},$ for any $t<1 $, we have
$$ \left( \wbrp_s , \lbrp_s , \hbrp_s
\right)_{0\leq s \leq t}
\xrightarrow[p\rightarrow \infty]{(d)}
\left( \xbr_s , L_s (\xbr ) , \hbr_s
\right)_{0\leq s\leq t }    .$$
\end{lemma}

\dem
Set $f(n,k)=\PP ( W_n =k) $ , $n\in \NN   ,  k\in \Z $. Let $F$ be any bounded continuous
functional on $\D ([0,t] , \R^3 )$.
The Markov property at time $[pt] $ under $\PP (.\mid W_p =-1 )$ implies that
\begin{equation}\begin{split}
\label{lemme11}
&\E \left[ F\left(
\wbrp_s , \lbrp_s , \hbrp_s
;  0\leq s\leq t \right) \right] \\
&\qquad=
\E \left[ \frac{ f( p-[pt] , -1-W_{[pt]} )}{f(p,-1)} \right. \\
&\qquad\qquad\quad\left. \times   F\left(
\frac{1}{a_p} W_{[ps]} , \frac{a_p}{p} L_{[ps]} (W), \frac{a_p}{p} H_{[ps]} (W)   ;
0\leq s\leq t \right) \right]   .
\end{split}\end{equation}

Since we assume \ref{attraction} and since $\mu $
(and thus $\nu $) is aperiodic, we can apply the Gnedenko local limit theorem
to $\nu $ in order to get
$$ \lim_{p\rightarrow \infty}\sup_{k\in \Z}\mid a_p
f(p-[pt] , k )- p_{1-t} (k/a_p)\mid=0$$
(see Bingham, Goldies and Teugels \cite{BiGoTe}). This result combined with (\ref{lemme11}),
the continuity of $x\rightarrow p_{1-t} (x)$ and Lemma
\ref{jointconv} gives
$$ \lim_{p\rightarrow \infty }
\E \left[ F\left(\wbrp_s , \lbrp_s , \hbrp_s;  0\leq s\leq t \right) \right] =
\E\left[ \frac{p_{1-t} (-X_t )}{p_1 (0)}
F \left( X_s , L_s , H_s ;  0\leq s\leq t\right)\right] $$
\noi and the lemma follows from (\ref{abscont}).\cqfd

Next, we need to prove the following lemma:
\begin{lemma}
\label{lemme2}
Under the assumptions of Theorem \ref{maintheo}, we have
$$ \left( \rlbrp , \rwbrp , \wbrp \right) \xrightarrow[p \rightarrow \infty ]{(d)}
\left( L(\rxbr ) , \rxbr , \xbr \right)
  .$$
\end{lemma}

\dem
First, let us show that $\wbrp $ converges to $\xbr $ in distribution in $ \D ([0,1], \R )$.
From Lemma \ref{lemme1} and the usual
tightness criterion, we only need to prove
\begin{equation}
\label{lemme21}
\lim_{\delta \rightarrow 0} \limsup_{p\rightarrow +\infty } \PP \left(
\sup_{s\in [1-\delta , 1 ]} \mid  \wbrp_s - \wbrp_1 \mid  > \eta \right) = 0
\end{equation}
\noi for any $\eta >0 $. Notice that the two
variables
$$ \sup_{[p(1-\delta )]\leq k \leq p }
\mid W_k -W_p \mid  \quad {\rm and} \quad
\sup_{0\leq k \leq p-[p(1-\delta )] } \mid W_k \mid
$$
\noi have the same law
under $\PP(  .  \mid W_p =-1 ) $. Thus
\begin{equation}
\label{lemme22}
\PP \left( \sup_{ s\in [1-\delta , 1 ]} \left|  \wbrp_s - \wbrp_1 \right|  > \eta \right)   \leq
\PP \left( \sup_{s \in [0, \delta +1/p ]} \mid  \wbrp_s\mid  > \eta \right)   .
\end{equation}
But Lemma \ref{lemme1} implies, for any $\eta >0 $,
$$ \lim_{\delta \rightarrow 0} \limsup_{p\rightarrow +\infty } \PP
\left( \sup_{s \in [0, \delta +1/p ]} \mid  \wbrp_s\mid  > \eta \right)
= 0   .$$
\noi Then, (\ref{lemme21}) follows from (\ref{lemme22}).

We now prove
\begin{equation}
\label{lemme24}
\lbrp \xrightarrow[p\rightarrow \infty]{(d)} L( \xbr )   .
\end{equation}
First, from Lemma \ref{lemme1} we have, for any $t<1 $,
$$(\lbrp_s )_{0\leq s \leq t}
\xrightarrow[p\rightarrow \infty]{(d)}
(L_s (\xbr ))_{0\leq s\leq t} $$
\noi in distribution
in $\D ([0,t] , \R)$. Next, recall
that $\PP$-a.s. there exists a small interval $(1-\delta , 1]$ on which
$L(\xbr ) $ is constant and equal to $L_1 (\xbr )$. So, we only need to prove
$$ \lim_{\delta \rightarrow 0} \limsup_{p\rightarrow +\infty } \PP \left(
\sup_{s\in [1-\delta , 1 ]} \left|  \lbrp_s - \lbrp_1 \right|> \eta \right) = 0 $$
\noi for any $\eta >0 $. But this is immediate from the observation that
$$  \lim_{\delta \rightarrow 0} \limsup_{p\rightarrow +\infty } \PP \left(
\sup_{s\in [1-\delta , 1 ]} W_s^{br , p} =
\sup_{s\in [0 , 1 ]} W_s^{br , p} \right) =0   .$$
which itself follows from the convergence of $ W^{br , p}$ to $X^{br }$.

Since $\wbrp $ (reps. $\xbr $) has the same law as $\rwbrp $ (resp. $\rxbr$),
the lemma is equivalent to
$$ \left( \lbrp , \wbrp , \rwbrp \right) \xrightarrow[p \rightarrow \infty ]{(d)}
\left( L(\xbr ) , \xbr , \rxbr \right)   .$$
First notice that the laws of $(\lbrp ,\wbrp , \rwbrp )$ are tight in the space of
all probability measures on
$\D([0,1], \R^3 )$. We only need to prove the convergence of the finite dimensional marginals.
By Lemma \ref{lemme1}, we see that the only
possible weak limit of the laws of $(\lbrp , \wbrp )$ is the law of $ (L(\xbr ) , \xbr ) $.
Since $\xbr $ has no fixed discontinuities,
we have for any $t_1, \ldots , t_n $ in  $[0,1]$
$$ \left( \lbrp_{t_i} , \wbrp_{t_i} \right)_{1\leq i\leq n} \xrightarrow[p \rightarrow \infty ]{(d)}
\left( L_{t_i}( \xbr) , \xbr_{t_i} \right)_{1\leq i\leq n}   . $$
For the same reason $\rxbr_t = -\xbr_{1-t} $, $\PP$-a.s. for any $t$ in $[0,1]$. So, we get
$$ \left( \lbrp_{t_i} , \wbrp_{t_i} ,\wbrp_{1} -\wbrp_{1-t_i}\right)_{1\leq i\leq n}
\xrightarrow[p \rightarrow \infty ]{(d)}
\left( L_{t_i}( \xbr) , \xbr_{t_i} , \rxbr_{t_i} \right)_{1\leq i\leq n}   . $$
But we have for any $t$ in $[0,1]$ the convergence in probability
$$ \wbrp_{1} -\wbrp_{1-t} -\rwbrp_t \xrightarrow[p \rightarrow \infty ]{ } 0 $$
because $\wbrp_{1} -\wbrp_{1-t} -\rwbrp_t $ has the same law as $(W_{p-[pt]} -W_{[p(1-t)]} )/a_p $
under $\PP (\cdot | W_p =-1 )$.
Thus, we have
$$ \left( \lbrp_{t_i} , \wbrp_{t_i} , \rwbrp_{t_i} \right)_{1\leq i\leq n} \xrightarrow[p
\rightarrow \infty ]{(d)}
\left( L_{t_i}( \xbr) , \xbr_{t_i} , \rxbr_{t_i} \right)_{1\leq i\leq n}   , $$
that implies the desired result.\cqfd

Next, we claim that the two following lemmas imply Proposition \ref{propaux}.
\begin{lemma}
\label{lemme3}
Under the assumptions of Theorem $\ref{maintheo},$ we have
$$ (\wbrp , \mbrp ) \xrightarrow[p\rightarrow \infty]{(d)}
 (\xbr , \mbr ) $$
\end{lemma}

\begin{lemma}
\label{lemme4}
Under the assumptions of Theorem $\ref{maintheo},$ the laws of the processes $( \hbrp )$ are tight
in the space of all probability measures on $\D( [0,1] , \R)$.
\end{lemma}

\bigskip

{\sc End of the proof of Proposition \ref{propaux}.}
The previous two lemmas imply that
the laws of $(\wbrp ,  \hbrp , \mbrp)$  are tight in the space of all probability measures
on $\D([0,1], \R^2 )$.
Let us assume that a subsequence of the sequence
$(\wbrp , \hbrp )$ converges in distribution
in $\D(\R_+ , \R^2 )$ to a certain process $(A, B)$. By Lemma \ref{lemme1}, it follows that
$$ (A_s , B_s )_{0\leq s\leq t}
\overset{\rm(law)}{=}
(\xbr_s , \hbr_s )_{0\leq s\leq t}  ,$$
for any $t<1 $. Also Lemma \ref{lemme3} implies $A \overset{\rm(law)}{=} X^{br }$.
Then, observe
that $\hbrp_1 = \rlbrp_1 $ , $p\geq 1 $ and that $\hbr_1 = L_1 (\rxbr )$.
From Lemma \ref{lemme2}, we get

$$ \left( A, B_1 \right) \overset{\rm(law)}{=} \left(  X^{br }, L_1 (\widehat{X}^{br } )
\right) = \left( X^{br } , H^{br }_1 \right)   .$$
This is more than enough to conclude that
$$ \left( A, B \right) \overset{\rm(law)}{=} \left( X^{br } , H^{br } \right)   .$$
So we have
$$ (\wbrp , \hbrp )  \xrightarrow[p\rightarrow \infty]{(d)}
(\xbr , \hbr )  .$$
Together with Lemma \ref{lemme3}, this implies that the only
possible weak limit of the laws of $(\wbrp ,  \hbrp , \mbrp)$ is the law of
$ (X^{br } , H^{br } , M^{br} )$. That completes the proof of Proposition \ref{propaux}.\cqfd

 {\sc Proof of Lemma \ref{lemme3} } We can apply Skorokhod's representation theorem
to replace the weak convergence of Lemma \ref{lemme2} by an a.s. convergence.
For convenience, we keep the same notation for the processes and the underlying probability
space, so we can suppose
\begin{equation}
\label{lemme31}
 (\rlbrp ,  \rwbrp , \wbrp )  \xrightarrow[p\rightarrow \infty]{\quad} ( L(\rxbr ) , \rxbr , \xbr )
\end{equation}
\noi $\PP$-a.s. for the Skorokhod topology in $\D( [0,1] , \R^3 )$.

For any $p\geq 1 $, we define the process
$(B^p_x ;   x\geq 0 )$ by
$$ B^p_x  = \rlbrp_1 -\rlbrp_{1\wedge  T_x (\rwbrp )}   .$$
\noi We get from the definition of $M$ the following inequality:
\begin{equation}
\label{inegal}
\sup_{t\in [0,1]}  \left| \mbrp_t -B^p_{-\inf_{[0,t]} \wbrp } \right| \leq    \app   .
\end{equation}
because
$$ \left| p\wedge T_{-\inf_{[0,t]} \wbrp } (\rwbrp ) - \gamma_p ( [pt] ) \right| \leq 1   .$$
We claim next that
\begin{equation}
\label{lemme32}
(B^p_x)_{x\geq 0 } \xrightarrow[p\rightarrow \infty]{\quad} (B_x )_{x\geq 0}   ,
\end{equation}
\noi $\PP$-a.s. for the Skorokhod topology in $\D (\R_+ , \R )$ :
Lemma 2.10, page 304, Chapter VI in \cite{Ja} shows for any $x \geq 0 $ that the functional
$1\wedge  T_x (.) $, defined on $\D ([0,1], \R)$ is continuous with respect to the Skorokhod
topology at any path $\omega $ satisfying
$ x\not\in J(\omega )$, where
$$ J (\omega ) = \{ y >0   :  T_{y+} (\omega )
 > T_y (\omega ) \}   .$$
An elementary argument shows that for any $x\geq 0$, $\PP ( x \in J(X) ) =0$. Then, we can
use the absolute continuity relation
(\ref{abscont}) to deduce  that $\PP$-a.s. $x$ is not in $J(\xbr ) $. Hence
\begin{equation}
\label{lemma33}
\PP\mbox{-a.s.,} \qquad 1 \wedge  T_q ( \rwbrp ) \longrightarrow 1\wedge  T_q (\rxbr ),\qquad  q\in \Q_+.
\end{equation}
\noi  Since $L(\rxbr )$ is continuous, a standard argument implies that $\PP$-a.s. $\rlbrp $
 converges to $ L(\rxbr )$ uniformly on $[0,1]$. Then, by (\ref{lemma33}), it follows that
\begin{equation}
\label{dense}
\left( B^p_{q_1 } , \ldots , B^p_{q_n } \right) \longrightarrow \left( B_{q_1 } ,
\ldots , B_{q_n } \right), \qquad \PP\mbox{-a.s.}
\end{equation}
\noi for any positive rational numbers $q_1 , q_2 , \ldots , q_n $. Next, observe that $B^p $ and $B$ are
nondecreasing processes and that $B$ is continuous [cf. Proposition
\ref{vervaat} (i))] so (\ref{dense}) implies the desired claim by a standard argument.

  It remains to prove that (\ref{lemme32}) implies the lemma:
Since $I^{br} $ is continuous, (\ref{lemme31})
implies that
$$ \lim_{p\rightarrow \infty} \sup_{t\in [0,1]}\left|
\inf_{[0, t]}  \wbrp  -
I^{br}_t  \right| =0   \qquad {\rm a.s.}  $$
\noi This, combined with (\ref{lemme32}), shows that $\PP$-a.s. $\mbrp $
converges to $\mbr$, uniformly on $[0,1]$. As $\mbr $ is continuous, a standard argument
(see \cite{Ja}, Proposition 1.23, page 293) implies that  $(\wbrp , \mbrp )$
converges almost surely to $(\xbr , \mbr )$  for the Skorokhod topology in $\D ([0,1] , \R^2 ) $.
That completes the proof of the lemma.\cqfd

\noi {\sc Proof of Lemma \ref{lemme4} }  By Lemma \ref{lemme1}, it is sufficient to show
\begin{equation}
\label{lemme41}
\lim_{\delta \rightarrow 0} \limsup_{p\rightarrow +\infty } \PP \left(
\sup_{s\in [1-\delta , 1 ]} \mid  \hbrp_s - \hbrp_1 \mid  > \eta \right) = 0
\end{equation}
\noi for any $\eta >0 $: Recall that
$$ G_p = \inf \left\{ 0\leq k\leq p   :  W_k= \inf_{0\leq j\leq k} W_j \right\} $$
and that
$$ W^{([p\varepsilon ])}_k = W_{[p\varepsilon ] +k} -W_{[p\varepsilon ]}   , \qquad k\geq 0   .$$
Let $\delta , \varepsilon >0 $ and set $A_p =\{
[p\varepsilon ] \leq G_p \leq p-[p\delta ] \} $. Observe that on $A_p $
$$\inf_{[p\varepsilon ] \leq i\leq [p\varepsilon ] +k } H_i (W) =0   ,\qquad
1-[p\delta ]-[ p \varepsilon ] \leq k \leq 1-[p\varepsilon ]   . $$
Then,  by (\ref{combi1}), it follows that
\begin{equation}
\label{nullite}
H_{ [p\varepsilon ] +k }(W) = H_k (W^{([p\varepsilon ] )} )   ,\quad
1-[p\delta ]- [p\varepsilon ]\leq k \leq 1-[p\varepsilon ]   .
\end{equation}
But it is easy to see that under $\PP ( \cdot \mid W_p=-1 )  $,
$ (W_i^{([p\varepsilon] )};  0\leq i\leq  1-[p\varepsilon ]) $ and
$ (W_i ;  0\leq i\leq  1-[p\varepsilon]) $ have the same law. Thus, by
(\ref{nullite}), we have, for any $a>0,$
\begin{equation}\begin{split}
\label{inegalite}
&\PP \left(
\sup_{1-[p\delta ]-[p\varepsilon ] \leq k \leq 1-[p\varepsilon ] }
\left| H_{ [p\varepsilon] +k }(W) - H_p (W) \right| >a   \arrowvert   W_p =-1 \right)\\
&\qquad    \leq   \PP\left( A^c_p \mid W_p =-1 \right) \\
&\qquad\quad+\PP \left( \sup_{1-[p\delta ]- [p\varepsilon ]\leq k \leq 1-[p\varepsilon ] }
\mid H_{ k }(W) - H_{p- [p\varepsilon ]} (W) \mid >a   \arrowvert
W_p =-1 \right)
\end{split}\end{equation}
Recall that under $\PP ( \cdot   \mid W_p =-1 ) $, the instant
$G_p $ is uniformly distributed on $\{ 1, \ldots , p\} $. So,
\begin{equation}
\label{nullite1}
\PP\left( A^c_p \mid W_p =-1 \right) \leq \delta +\varepsilon   .
\end{equation}
Set $\delta_p =[p\delta ]/p $ and $\varepsilon_p =[p\varepsilon ]/p $, and take $a= a_p\eta / p $ in
(\ref{inegalite}) in order to get
\begin{equation}\begin{split}
\label{inegalite2}
&\PP \left( \sup_{s\in [1-\delta , 1 ]} \mid  \hbrp_s - \hbrp_1 \mid  > \eta \right) \\
&\qquad\leq \PP \left( \sup_{s\in [1-\varepsilon_p -\delta_p , 1-\varepsilon_p ]} \mid  \hbrp_s -
\hbrp_{1-\varepsilon_p } \mid  > \eta \right) + \delta + \varepsilon.
\end{split}\end{equation}
\noi By Lemma \ref{lemme1}, it follows that
$$ \lim_{\delta \rightarrow 0} \limsup_{p\rightarrow +\infty }
\PP \left( \sup_{s\in [1-\varepsilon_p -\delta_p , 1-\varepsilon_p ]} \mid  \hbrp_s -
\hbrp_{1-\varepsilon_p } \mid >\eta \right) = 0.$$
Thus
$$  \limsup_{\delta \rightarrow 0} \limsup_{p\rightarrow +\infty }
\PP \left( \sup_{s\in [1-\delta,1]} \mid  \hbrp_s - \hbrp_1 \mid  > \eta \right) \leq
\varepsilon   ,$$
which yields the desired result by letting $\varepsilon $ go to $0$.\cqfd

\end{document}